\documentclass[aps,prb]{revtex4-2}

\usepackage{amsmath,amssymb,amsthm}
\usepackage{mathtools}
\usepackage{xspace}
\usepackage{braket}
\usepackage{hyphenat}
\usepackage{bbold}
\usepackage{enumerate}
\usepackage{tabularx}

\usepackage[british]{babel}

\usepackage{srcltx}
\usepackage{color}

\usepackage{graphicx}
\usepackage{amsmath}
\usepackage{amssymb}
\usepackage{bm}

\begin{document}

\title{On an infinite number of nonlinear Euler sums}
\author{J. Braun}
\affiliation{Ludwig Maximilians-Universit{\"a}t, M{\"u}nchen, Germany}
\author{H. J. Bentz}
\affiliation{Institut f\"{u}r Mathematik und Informatik, Universit\"{a}t Hildesheim, Germany}
\date{\today}

\begin{abstract}
Linear harmonic number sums had been studied by a variety of authors during the last centuries, but only few results are known about nonlinear Euler sums of quadratic
or even higher degree. The first systematic study on nonlinear Euler sums consisting of products of hyperharmonic sums had been published by Flajolet and Salvy in
1997 followed by similar studies presented during the last years by different authors. Although these studies had been restricted to sums where the nominator
consists of a product of even or odd hyperharmonic sums, where the denominator is of the type $1/k^n$. We have generalized these results to nonlinear Euler sums with
different denominators and nominators which consist in addition of mixed products between even and odd hyperharmonic numbers. In detail we present eight families of
quadratic Euler sums which are expressible by zeta values and special types of linear Euler sums only where the order of the nonlinear Euler sums is always an odd 
number. The resulting eight different families of nonlinear Euler sums which we discovered consist of various products between even and odd hyperharmonic numbers,
divided by three different types of denominators $1/k^n$, $1/((2k-1)^n)$ and $1/(k(2k-1))$. The calculational scheme is based on proper two-valued integer functions, which
allow us to compute these sequences explicitly in terms of zeta values and pairs of odd-type linear harmonic numbers and even hyperharmonic numbers of second order.
\end{abstract}
\maketitle
%Keywords: Harmonic Numbers, Nonlinear Euler sums, Riemann Zeta function

\section{Introduction}
Harmonic numbers and their generalizations called hyper-harmonic numbers are defined by
\begin{eqnarray}
H_k = H_k^{(1)} = \sum^{k}_{i=1} \frac{1}{i}~, ~~~~~~~~~~~~  H_k^{(n)} = \sum^{k}_{i=1} \frac{1}{i^n}~.
\end{eqnarray}
It has been discovered by many authors in the past starting with Goldbach 1742 and later on with Euler, Ramanujan, Borwein, Salvi and others that linear sums of the type
\begin{eqnarray}
t(n,m) = \sum^{\infty}_{k=1} \frac{H_k^{(n)}}{k^m} 
\end{eqnarray}
are expressible in terms of zeta values only for an odd weight $p=n+m$, where $p=n+m$ defines the order of the corresponding sum. An excellent overview about these works
is found in the paper of Flajolet and Salvy \cite{fla97} and in citations therein. Further work on similar types of Euler sums can be found, for example, in
\cite{ce16,ce17,ce18,ce20}. Although, only few sums of the type
\begin{eqnarray}
s(n,m) = \sum^{\infty}_{k=1} \frac{h_k^{(n)}}{k^m} 
\end{eqnarray}
had been discovered in the past \cite{sit85,zeh07,ade16,ade16a,ce20,bra21}. In analogy these sums are expressible in terms of zeta values only for an odd weight $p=n+m$,
where again $p=n+m$ defines the order of the corresponding sum. Here the symbol $h_k$ denote odd-type (hyper-)harmonic sums which are defined as follows:
\begin{eqnarray}
h_k = h_k^{(1)} = \sum^{k}_{i=1} \frac{1}{2i-1}~,~~~~~~~~~~~~  h_k^{(n)} = \sum^{k}_{i=1} \frac{1}{(2i-1)^n}~.
\end{eqnarray} 
This way the subject of the present paper is to present a formalism that allows for an explicit calculation of quadratic Euler sums of odd order in terms zeta values and
Polylogarithms. In detail we introduce a calculational scheme which is based on proper two-valued integer functions \cite{bra20,bra21,bra22}. In this sense our work can be
seen as a generalization of the paper that we published in 2022 for even order Euler sums \cite{bra22}. We will show that the corresponding Euler sums can be expressed
in terms of zeta values and Polylogarithmic values $Li_n(1/2)$.

\section{The order three case}
In this section we will show that the following families of quadratic Euler sums can be calculated in terms of zeta values only.  
We start with the limiting cases, which means with the corresponding linear Euler sums. It follows:
\begin{eqnarray}
\sum^{\infty}_{k=1} \frac{H^{(2)}_k}{k(2k-1)} = -\frac{5}{2}\zeta(3) + 2ln(2)\zeta(2) - 2\zeta(2) + 8ln(2)~.
\end{eqnarray}
Analogously it follows:
\begin{eqnarray}
\sum^{\infty}_{k=1} \frac{h^{(2)}_k}{k(2k-1)} = \frac{21}{8}\zeta(3) - \frac{3}{2}ln(2)\zeta(2)~.
\end{eqnarray}

\subsection{Third family}
The following type of nonlinear Euler sums of order 3 can be expressed in terms of zeta values only:
\begin{eqnarray}
\sum^{\infty}_{k=1} \frac{H^{(a)}_k H^{(b)}_k}{k(2k-1)}
\end{eqnarray}
with $a,b \in \mathbb{N}$, $a+b=1$. We get: 
\begin{eqnarray}
\sum^{\infty}_{k=1} \frac{H_k H_k}{k(2k-1)} = -\frac{1}{2}\zeta(3) - 2ln(2)\zeta(2) + 2\zeta(2) + 8ln(2) - 8\left(ln(2)\right)^2 +
\frac{8}{3}\left(ln(2)\right)^3~,
\end{eqnarray}
where we used the identity (Eq.~(84)) from \cite{bra22}.

\subsection{Fifth family}
The following type of nonlinear Euler sums of order 3 can be expressed in terms of zeta values only:
\begin{eqnarray}
\sum^{\infty}_{k=1} \frac{H^{(a)}_k h^{(b)}_k}{k(2k-1)}
\end{eqnarray}
with $a,b \in \mathbb{N}$, $a+b=1$. We get: 
\begin{eqnarray}
\sum^{\infty}_{k=1} \frac{H_k h_k}{k(2k-1)} = \frac{5}{4}\zeta(3) - 2ln(2)\zeta(2) + 2\zeta(2)~.
\end{eqnarray}
where we used the identity
\begin{eqnarray}
\sum^{\infty}_{k=1} \frac{1}{k(2k-1)} \left( \sum^k_{i=1} \frac{h_i}{i} \right) = 2\zeta(2) - \zeta(3)~.
\end{eqnarray}

\subsection{Sixth family}
The following type of nonlinear Euler sums of order 3 can be expressed in terms of zeta values only:
\begin{eqnarray}
\sum^{\infty}_{k=1} \frac{h^{(a)}_k h^{(b)}_k}{k(2k-1)}
\end{eqnarray}
with $a,b \in \mathbb{N}$, $a+b=1$. We get: 
\begin{eqnarray}
\sum^{\infty}_{k=1} \frac{h_k h_k}{k(2k-1)} = \frac{3}{8}\zeta(3) + \frac{3}{2}ln(2)\zeta(2),
\end{eqnarray}
where we used lemma 4a from \cite{bra22}. 

\section{The order five case}
In the following section we demonstrate that the different families of quadratic Euler sums can be calculated in terms of zeta values and/or Polylogarithmic values.  

\subsection{First family}
The following type of nonlinear Euler sums of order 5 can be expressed in terms of zeta values and the Polylogarithmic values $Li_4(1/2)$ and $Li_5(1/2)$:
\begin{eqnarray}
\sum^{\infty}_{k=1} \frac{H^{(a)}_k h^{(b)}_k}{k^c}
\end{eqnarray}
with $a,b,c \in \mathbb{N}$, $a+b+c=$odd. First we get: 
\begin{eqnarray}
\sum^{\infty}_{k=1} \frac{H_k h_k}{k^3} = \frac{279}{16}\zeta(5) - 7\left(ln(2)\right)^2\zeta(3) + \frac{8}{3}\left(ln(2)\right)^3\zeta(2) - 
\frac{8}{15} \left(ln(2)\right)^5 - 16ln(2)Li_4\left(\frac{1}{2}\right) - 16 Li_5\left(\frac{1}{2}\right)~.
\end{eqnarray}
Furthermore we get:
\begin{eqnarray}
\sum^{\infty}_{k=1} \frac{H^{(2)}_k h_k}{k^2} = \frac{93}{16}\zeta(5) - \frac{7}{4}\zeta(2)\zeta(3)~,
\end{eqnarray}
and finally
\begin{eqnarray}
\sum^{\infty}_{k=1} \frac{H_k h^{(2)}_k}{k^2} &=& -\frac{589}{32}\zeta(5) + \frac{21}{8}\zeta(2)\zeta(3) + 7\left(ln(2)\right)^2\zeta(3) 
- \frac{8}{3}\left(ln(2)\right)^3\zeta(2) + \frac{8}{15} \left(ln(2)\right)^5 + 16ln(2)Li_4\left(\frac{1}{2}\right) \nonumber \\ &+& 
16 Li_5\left(\frac{1}{2}\right)~.
\end{eqnarray}
With this all members of the first family of order five have been explicitly calculated.

\subsection{Second family}
The following type of nonlinear Euler sums of order 5 can be expressed in terms of zeta values and the Polylogarithmic values $Li_4(1/2)$ and $Li_5(1/2)$:
\begin{eqnarray}
\sum^{\infty}_{k=1} \frac{H^{(a)}_k H^{(b)}_k}{(2k-1)^c}
\end{eqnarray}
with $a,b,c \in \mathbb{N}$, $a+b+c=$odd. First we get: 
\begin{eqnarray}
\sum^{\infty}_{k=1} \frac{H_k H_k}{(2k-1)^3} &=& \frac{31}{8}\zeta(5) - \frac{7}{8} \zeta(2)\zeta(3) - \frac{45}{8}ln(2)\zeta(4) + \frac{45}{8}\zeta(4) + 
\frac{7}{2}\left(ln(2)\right)^2\zeta(3) - 7ln(2) \zeta(3) - \frac{7}{2}\zeta(3) + 6ln(2)\zeta(2) \nonumber \\ &-& 5\zeta(2) + 12 ln(2) -4\left(ln(2)\right)^2~.
\end{eqnarray}
This Euler sum has been introduced first in \cite{bra22}.
The calculation for the second sum of oder five gives the result:  
\begin{eqnarray}
\sum^{\infty}_{k=1} \frac{H_k H^{(2)}_k}{(2k-1)^2} &=& -\frac{589}{32}\zeta(5) + \frac{121}{8}ln(2)\zeta(4) - \frac{121}{8}\zeta(4) + \frac{7}{4} \zeta(2)\zeta(3)
- 7\left(ln(2)\right)^2\zeta(3) + 14ln(2)\zeta(3) + \frac{21}{2}\zeta(3) \nonumber \\ &+& \frac{4}{3}\left(ln(2)\right)^3\zeta(2) - 4\left(ln(2)\right)^2\zeta(2) 
- 8ln(2)\zeta(2) + 6\zeta(2) -24ln(2) + 8\left(ln(2)\right)^2 + \frac{2}{3}\left(ln(2)\right)^4 \nonumber \\ &-& \frac{2}{15}\left(ln(2)\right)^5 +
16Li_4\left(\frac{1}{2}\right) + 16Li_5\left(\frac{1}{2}\right)~.
\end{eqnarray}
With this all members of the second family of order five have been explicitly calculated.

\subsection{Third family}
The following type of nonlinear Euler sums of order 5 can be expressed in terms of zeta values and the Polylogarithmic values $Li_4(1/2)$ and $Li_5(1/2)$:
\begin{eqnarray}
\sum^{\infty}_{k=1} \frac{H^{(a)}_k H^{(b)}_k}{k(2k-1)}
\end{eqnarray}
with $a,b,c \in \mathbb{N}$, $a+b=$even. First we get: 
\begin{eqnarray}
\sum^{\infty}_{k=1} \frac{H_k H^{(3)}_k}{k(2k-1)} &=& \frac{251}{8}\zeta(5) - \frac{53}{2}ln(2)\zeta(4) + \frac{53}{2}\zeta(4) - \frac{3}{2} \zeta(2)\zeta(3) + 
12\left(ln(2)\right)^2\zeta(3) - 24ln(2)\zeta(3) - 8\zeta(3) \nonumber \\ &-& \frac{8}{3}\left(ln(2)\right)^3\zeta(2) + 
8\left(ln(2)\right)^2\zeta(2) + 32 ln(2) - 16\left(ln(2)\right)^2 - \frac{4}{3}\left(ln(2)\right)^4 + \frac{4}{15}\left(ln(2)\right)^5 \nonumber \\ &-&~
32Li_4\left(\frac{1}{2}\right) - 32Li_5\left(\frac{1}{2}\right),
\end{eqnarray}
and for the second sum it follows: 
\begin{eqnarray}
\sum^{\infty}_{k=1} \frac{H^{(2)}_k H^{(2)}_k}{k(2k-1)} &=& \frac{5}{4}\zeta(5) + 5ln(2)\zeta(4) - 5\zeta(4) - 3\zeta(2)\zeta(3) - 16\zeta(3) + 
16ln(2)\zeta(2) - 8\zeta(2) + 32ln(2)~.
\end{eqnarray}
With this all members of the third family of order five have been explicitly calculated.

\subsection{Fourth family}
The following type of nonlinear Euler sums of order 5 can be expressed in terms of zeta values and the Polylogarithmic values $Li_4(1/2)$ and $Li_5(1/2)$:
\begin{eqnarray}
\sum^{\infty}_{k=1} \frac{h^{(a)}_k h^{(b)}_k}{(2k-1)^c}
\end{eqnarray}
with $a,b,c \in \mathbb{N}$, $a+b+c=$odd. First we get: 
\begin{eqnarray}
\sum^{\infty}_{k=1} \frac{h_k h_k}{(2k-1)^3} &=& -\frac{217}{128}\zeta(5) + \frac{83}{32}ln(2)\zeta(4) - \frac{3}{32}\zeta(2)\zeta(3) + \frac{1}{6}\left(ln(2)\right)^3\zeta(2)
-\frac{1}{60}\left(ln(2)\right)^5 + 2Li_5\left(\frac{1}{2}\right)~.
\end{eqnarray}
and for the second sum it follows: 
\begin{eqnarray}
\sum^{\infty}_{k=1} \frac{h_k h^{(2)}_k}{(2k-1)^2} &=& \frac{93}{256}\zeta(5) + \frac{75}{64}ln(2)\zeta(4) + \frac{3}{32}\zeta(2)\zeta(3)~.
\end{eqnarray}
With this all members of the fourth family of order five have been explicitly calculated.

\subsection{Fifth family}
The following type of nonlinear Euler sums of order 5 can be expressed in terms of zeta values and the Polylogarithmic values $Li_4(1/2)$ and $Li_5(1/2)$:
\begin{eqnarray}
\sum^{\infty}_{k=1} \frac{H^{(a)}_k h^{(b)}_k}{(2k-1)^c}
\end{eqnarray}
with $a,b,c \in \mathbb{N}$, $a+b+c=$odd. First we get: 
\begin{eqnarray}
\sum^{\infty}_{k=1} \frac{H^{(2)}_k h_k}{(2k-1)^2} &=& \frac{279}{32}\zeta(5) - \frac{121}{16}ln(2)\zeta(4) + \frac{63}{16}\zeta(2)\zeta(3) + \frac{7}{2}\zeta(3) +
\frac{2}{3}\left(ln(2)\right)^3\zeta(2) + 3ln(2)\zeta(2) - 4\zeta(2) \nonumber \\ &-& \frac{1}{5}\left(ln(2)\right)^5 - 8ln(2)Li_4\left(\frac{1}{2}\right) - 
16Li_5\left(\frac{1}{2}\right)~.
\end{eqnarray}
and for the second sum it follows: 
\begin{eqnarray}
\sum^{\infty}_{k=1} \frac{H_k h_k}{(2k-1)^3} &=& \frac{279}{128}\zeta(5) - \frac{19}{16}ln(2)\zeta(4) + \frac{83}{32}\zeta(4) - \frac{7}{8}\zeta(2)\zeta(3) -
\frac{7}{8}\zeta(3) - \frac{1}{2}\left(ln(2)\right)^3\zeta(2) + \frac{1}{2}\left(ln(2)\right)^2\zeta(2) \nonumber \\ &-& \frac{3}{2}ln(2)\zeta(2) + \zeta(2) - 
\frac{1}{12}\left(ln(2)\right)^4 + \frac{1}{12}\left(ln(2)\right)^5 - 2Li_4\left(\frac{1}{2}\right) + 2ln(2)Li_4\left(\frac{1}{2}\right)~.
\end{eqnarray}
The third sum has been calculated to:
\begin{eqnarray}
\sum^{\infty}_{k=1} \frac{H_k h^{(2)}_k}{(2k-1)^2} &=& \frac{217}{64}\zeta(5) - \frac{75}{32}ln(2)\zeta(4) + \frac{75}{32}\zeta(4) - \frac{21}{32}\zeta(2)\zeta(3) -
\frac{21}{8}\zeta(3) + \frac{3}{2}ln(2)\zeta(2)~.
\end{eqnarray}
With this all members of the fifth family of order five have been explicitly calculated.

\subsection{Sixth family}
The following type of nonlinear Euler sums of order 5 can be expressed in terms of zeta values and the Polylogarithmic values $Li_4(1/2)$ and $Li_5(1/2)$:
\begin{eqnarray}
\sum^{\infty}_{k=1} \frac{h^{(a)}_k h^{(b)}_k}{k(2k-1)}
\end{eqnarray}
with $a,b \in \mathbb{N}$, $a+b=$even. First we get: 
\begin{eqnarray}
\sum^{\infty}_{k=1} \frac{h_k h^{(3)}_k}{k(2k-1)} &=& \frac{341}{128}\zeta(5) - \frac{23}{32}ln(2)\zeta(4) + \frac{9}{32}\zeta(2)\zeta(3) -
\frac{1}{6}\left(ln(2)\right)^3\zeta(2) + \frac{1}{60}\left(ln(2)\right)^5 - 2Li_5\left(\frac{1}{2}\right)~,
\end{eqnarray}
and for the second sum it follows:
\begin{eqnarray}
\sum^{\infty}_{k=1} \frac{h^{(2)}_k h^{(2)}_k}{k(2k-1)} &=& \frac{155}{64}\zeta(5) - \frac{45}{16}ln(2)\zeta(4) + \frac{9}{16}\zeta(2)\zeta(3)~.
\end{eqnarray}
With this all members of the sixth family of order five have been explicitly calculated.

\subsection{Seventh family}
The following type of nonlinear Euler sums of order 5 can be expressed in terms of zeta values and the Polylogarithmic values $Li_4(1/2)$ and $Li_5(1/2)$:
\begin{eqnarray}
\sum^{\infty}_{k=1} \frac{H^{(a)}_k h^{(b)}_k}{k(2k-1)}
\end{eqnarray}
with $a,b \in \mathbb{N}$, $a+b=$even. First we get: 
\begin{eqnarray}
\sum^{\infty}_{k=1} \frac{H^{(3)}_k h_k}{k(2k-1)} &=& \frac{713}{16}\zeta(5) - \frac{53}{4}ln(2)\zeta(4) + \frac{53}{4}\zeta(4) - \frac{31}{4}\zeta(2)\zeta(3) -
14ln(2)\zeta(3) -7\zeta(3) + \frac{4}{3}\left(ln(2)\right)^3\zeta(2) \nonumber \\ &+& 4\left(ln(2)\right)^2\zeta(2) + 8\zeta(2) - 
\frac{2}{3}\left(ln(2)\right)^4 - \frac{2}{5}\left(ln(2)\right)^5 - 16ln(2)Li_4\left(\frac{1}{2}\right)- 16Li_4\left(\frac{1}{2}\right) \nonumber \\ &-&
32Li_5\left(\frac{1}{2}\right)~.
\end{eqnarray}
For the second sum it follows:
\begin{eqnarray}
\sum^{\infty}_{k=1} \frac{H^{(2)}_k h^{(2)}_k}{k(2k-1)} &=& -\frac{279}{4}\zeta(5) + 34ln(2)\zeta(4) - \frac{151}{8}\zeta(4) + \frac{7}{4}\zeta(2)\zeta(3) -
14ln(2)\zeta(3) + \frac{21}{2}\zeta(3) - \frac{8}{3}\left(ln(2)\right)^3\zeta(2) \nonumber \\ &-& 4\left(ln(2)\right)^2\zeta(2) - 6ln(2)\zeta(2) + 
\frac{2}{3}\left(ln(2)\right)^4 + \frac{4}{5}\left(ln(2)\right)^5 + 32ln(2)Li_4\left(\frac{1}{2}\right) + 16Li_4\left(\frac{1}{2}\right) \nonumber \\ &+&
64Li_5\left(\frac{1}{2}\right)~.
\end{eqnarray}
The third sum has been calculated to:
\begin{eqnarray}
\sum^{\infty}_{k=1} \frac{H_k h^{(3)}_k}{k(2k-1)} &=& \frac{31}{64}\zeta(5) - \frac{67}{16}ln(2)\zeta(4) + \frac{11}{8}\zeta(4) + \frac{21}{16}\zeta(2)\zeta(3) +
\left(ln(2)\right)^3\zeta(2) - \left(ln(2)\right)^2\zeta(2) \nonumber \\ &+& \frac{1}{6}\left(ln(2)\right)^4 - \frac{1}{6}\left(ln(2)\right)^5 - 
4ln(2)Li_4\left(\frac{1}{2}\right) + 4Li_4\left(\frac{1}{2}\right)~.
\end{eqnarray}
With this all members of the seventh family of order five have been explicitly calculated.

\subsection{Eighth family}
The following type of nonlinear Euler sums of order 5 can be expressed in terms of zeta values and the Polylogarithmic values $Li_4(1/2)$ and $Li_5(1/2)$:
\begin{eqnarray}
\sum^{\infty}_{k=1} \frac{h^{(a)}_k h^{(b)}_k}{k^c}
\end{eqnarray}
with $a,b,c \in \mathbb{N}$, $a+b+c=$odd. First we get: 
\begin{eqnarray}
\sum^{\infty}_{k=1} \frac{h_k h^{(2)}_k}{k^2} &=& \frac{31}{8}\zeta(5) - \frac{7}{8}\zeta(2)\zeta(3)~.
\end{eqnarray}
For the second sum it follows:
\begin{eqnarray}
\sum^{\infty}_{k=1} \frac{h_k h_k}{k^3} &=& -\frac{31}{16}\zeta(5) + \frac{7}{4}\zeta(2)\zeta(3)~.
\end{eqnarray}
This Euler sum has been calculated first in \cite{bra20}.

\section{The order seven case}
In the following section we demonstrate that the different families of quadratic Euler sums can be calculated in terms of zeta values and/or Polylogarithmic values.  

\subsection{First family}
The order seven case seems to be more complicated as higher order polylogarithmic values like  $Li_7(1/2)$ are not known explicitly but most probably necessary for 
the complete calculation of the nonlinear Euler sums being members of the first family. Nevertheless a significant reduction is possible. This will be shown in
the following.
\begin{eqnarray}
\sum^{\infty}_{k=1} \frac{H^{(a)}_k h^{(b)}_k}{k^c}
\end{eqnarray}
with $a,b,c \in \mathbb{N}$, $a+b+c=$odd. First we get: 
\begin{eqnarray}
\sum^{\infty}_{k=1} \frac{H_k^{(3)} h_k}{k^3} = \frac{5715}{64}\zeta(7) - \frac{279}{8}\zeta(2)\zeta(5) - \frac{329}{16}\zeta(3)\zeta(4) + 
\zeta(3)\sum^{\infty}_{k=1} \frac{h_k}{k^3}- 
\frac{7}{2} \sum^{\infty}_{k=1} \frac{H_k^{(2)} h_k}{k^4}.
\end{eqnarray}
and
\begin{eqnarray}
\sum^{\infty}_{k=1} \frac{H_k^{(4)} h_k}{k^2} = \frac{635}{16}\zeta(7) - \frac{31}{2}\zeta(2)\zeta(5) - \frac{63}{8}\zeta(3)\zeta(4) -
\sum^{\infty}_{k=1} \frac{H_k^{(2)} h_k}{k^4}.
\end{eqnarray}
and
\begin{eqnarray}
\sum^{\infty}_{k=1} \frac{H_k h_k^{(2)}}{k^4} = \frac{127}{32}\zeta(7) + \frac{31}{4}\zeta(2)\zeta(5) - \frac{41}{4}\zeta(3)\zeta(4) -
2\sum^{\infty}_{k=1} \frac{H_k h_k}{k^5} - \frac{1}{2}\sum^{\infty}_{k=1} \frac{H_k^{(2)} h_k}{k^4}.
\end{eqnarray}
and
\begin{eqnarray}
\sum^{\infty}_{k=1} \frac{H_k h_k^{(3)}}{k^3} = -\frac{3937}{256}\zeta(7) - \frac{341}{32}\zeta(2)\zeta(5) + \frac{1687}{64}\zeta(3)\zeta(4) -
\frac{7}{4}\zeta(3)\sum^{\infty}_{k=1} \frac{h_k}{k^3} + \frac{3}{2}\sum^{\infty}_{k=1} \frac{H_k h_k}{k^5} + \frac{13}{8}\sum^{\infty}_{k=1} \frac{H_k^{(2)} h_k}{k^4}.
\end{eqnarray}
and
\begin{eqnarray}
\sum^{\infty}_{k=1} \frac{H_k h_k^{(4)}}{k^2} = \frac{889}{128}\zeta(7) + \frac{279}{32}\zeta(2)\zeta(5) - \frac{1001}{64}\zeta(3)\zeta(4) +
\frac{7}{4}\zeta(3)\sum^{\infty}_{k=1} \frac{h_k}{k^3} - \frac{1}{2}\sum^{\infty}_{k=1} \frac{H_k h_k}{k^5} - \frac{9}{8}\sum^{\infty}_{k=1} \frac{H_k^{(2)} h_k}{k^4}.
\end{eqnarray}
and
\begin{eqnarray}
\sum^{\infty}_{k=1} \frac{H_k^{(2)} h_k^{(2)}}{k^3} = -\frac{3175}{32}\zeta(7) + \frac{341}{8}\zeta(2)\zeta(5) + \frac{81}{4}\zeta(3)\zeta(4) +
2\sum^{\infty}_{k=1} \frac{H_k^{(2)} h_k}{k^4}.
\end{eqnarray}
and
\begin{eqnarray}
\sum^{\infty}_{k=1} \frac{H_k^{(3)} h_k^{(2)}}{k^2} = -\frac{9271}{64}\zeta(7) + \frac{527}{8}\zeta(2)\zeta(5) + \frac{193}{4}\zeta(3)\zeta(4) -
\zeta(3)\sum^{\infty}_{k=1} \frac{h_k}{k^3} + 5\sum^{\infty}_{k=1} \frac{H_k^{(2)} h_k}{k^4}.
\end{eqnarray}
and finally
\begin{eqnarray}
\sum^{\infty}_{k=1} \frac{H_k^{(2)} h_k^{(3)}}{k^2} = \frac{16637}{128}\zeta(7) - 62\zeta(2)\zeta(5) - \frac{119}{8}\zeta(3)\zeta(4) -
\frac{7}{2}\sum^{\infty}_{k=1} \frac{H_k^{(2)} h_k}{k^4}.
\end{eqnarray}
The corresponding mathematical reduction procedure which is behind it appears very inhomogeneous as one needs a variety of double valued help functions. This way we restrict here on
some representative examples. One can observe that the sum of the two Euler functions appearing in Eq.~(43) is trivial in the sense that only zeta values appear. We start with
the following identity:
\begin{eqnarray}
h_k = \frac{1}{2}H_k + \sum^{k}_{i=1}\frac{1}{i+k}~.
\end{eqnarray}
From this we get:
\begin{eqnarray}
\sum^{\infty}_{k=1} \frac{H_k^{(2)} h_k}{k^4} = \frac{1}{2} \frac{H_k H_k^{(2)}}{k^4} + \sum^{\infty}_{k=1} \frac{H_k^{(2)}}{k^4} \left(\sum^{k}_{i=1}\frac{1}{i+k}\right)~.
\end{eqnarray}
Rearranging the summations order it follows:
\begin{eqnarray}
\sum^{\infty}_{k=1} \frac{H_k^{(2)}}{k^4} \left(\sum^{k}_{i=1}\frac{1}{i+k}\right) &=& \frac{1}{2}\sum^{\infty}_{k=1} \frac{H_k^{(2)}}{k^5} + \sum^{\infty}_{i=1} \left(
\frac{1}{i} \sum^{\infty}_{k=1} \frac{H_{k+i}^{(2)}}{(k+i)^{4}}\right) - \sum^{\infty}_{i=1} \left( \frac{1}{i^2} \sum^{\infty}_{k=1} \frac{H_{k+i}^{(2)}}{(k+i)^{3}}\right)
\nonumber \\ &+& \sum^{\infty}_{i=1} \left( \frac{1}{i^3} \sum^{\infty}_{k=1} \frac{H_{k+i}^{(2)}}{(k+i)^{2}}\right) - \sum^{\infty}_{i=1} \left( \frac{1}{i^3}
\sum^{\infty}_{k=1} \frac{H_{k+i}^{(2)}}{(k+i)(k+2i)}\right)~.
\end{eqnarray}
With
\begin{eqnarray}
\frac{H_{k+i}^{(2)}}{(k+i)^{n}} = \frac{H_k^{(2)}}{k^n} - \sum^{i}_{k=1} \frac{H_k^{(2)}}{k^n}
\end{eqnarray}
the computation gives:
\begin{eqnarray}
\sum^{\infty}_{k=1} \frac{H_k^{(2)} h_k}{k^4} = \frac{3}{2} \frac{H_k H_k^{(2)}}{k^4} - \frac{1}{2} \frac{H_k^{(2)}}{k^5} - \frac{H_k^{(2)} H_k^{(2)}}{k^3} +
\frac{H_k^{(2)} H_k^{(3)}}{k^2} - \sum^{\infty}_{i=1} \left( \frac{1}{i^3} \sum^{\infty}_{k=1} \frac{H_{k+i}^{(2)}}{(k+i)(k+2i)}\right)~.
\end{eqnarray}
It remains to calculate
\begin{eqnarray}
\sum^{\infty}_{k=1} \frac{H_{k+i}^{(2)}}{(k+i)(k+2i)} = \sum^{\infty}_{k=1} \frac{H_k^{(2)}}{k(i+k)} - \sum^{i}_{k=1}\frac{H_k^{(2)}}{k(i+k)} 
\end{eqnarray}
explicitly. We found:
\begin{eqnarray}
\sum^{\infty}_{k=1} \frac{H_k^{(2)}}{k(i+k)} = \zeta(3)\frac{1}{i} - \zeta(2)\frac{1}{i^2} + \zeta(2)\frac{H_i}{i} - \frac{1}{i} \sum^{i-1}_{k=1}\frac{H_k}{k^2}
\end{eqnarray}
and
\begin{eqnarray}
\sum^{i}_{k=1}\frac{H_k^{(2)}}{k(i+k)} = 2\frac{H_i}{i} + \frac{H_i}{i^2} + 2\frac{H_i^{(3)}}{i} - 4\frac{h_i}{i} + \frac{3}{2}\frac{H_i H_i^{(2)}}{i} - \frac{H_i^{(2)} h_i}{i}
+ \frac{2}{i}\sum^{i}_{k=1} \frac{h_k}{k^2} - \frac{3}{i}\sum^{i}_{k=1} \frac{H_k}{k^2}~.
\end{eqnarray}
Finally using the two identities one gets:
\begin{eqnarray}
\sum^{\infty}_{k=1} \frac{H_k^{(2)} h_k}{k^4} &=& \zeta(2)\zeta(5) - \zeta(3)\zeta(4) - \frac{1}{2}\sum^{\infty}_{k=1} \frac{H_k^{(2)}}{k^5} - \zeta(2)\sum^{\infty}_{k=1}\frac{H_k}{k^4}
+ \zeta(4)\sum^{\infty}_{k=1}\frac{H_k}{k^2} + 2\sum^{\infty}_{k=1}\frac{H_k}{k^4} + \sum^{\infty}_{k=1}\frac{H_k}{k^5} - 4\sum^{\infty}_{k=1}\frac{h_k}{k^4}  \nonumber \\ &+& 
2\frac{H_k^{(3)}}{k^4} - \frac{H_k H_k^{(4)}}{k^2} + 3\frac{H_k H_k^{(2)}}{k^4} - \frac{H_k^{(2)} H_k^{(2)}}{k^3} + \frac{H_k^{(2)} H_k^{(3)}}{k^2} - 3\sum^{\infty}_{k=1}
\left(\frac{1}{k^4}\sum^{k}_{i=1}\frac{H_i}{i^2}\right) + 2\sum^{\infty}_{k=1}\left(\frac{1}{k^4}\sum^{k}_{i=1}\frac{h_i}{i^2}\right)  \nonumber \\ &-&
\sum^{\infty}_{k=1} \frac{H_k^{(2)} h_k}{k^4}~.
\end{eqnarray}
Using 
\begin{eqnarray}
\sum^{\infty}_{k=1} \left(\frac{1}{k^4}\sum^{k}_{i=1}\frac{H_i}{i^2}\right) = \zeta(4)\sum^{\infty}_{k=1}\frac{H_k}{k^2} + \sum^{\infty}_{k=1}\frac{H_k}{k^6} -
\sum^{\infty}_{k=1} \frac{H_k H_k^{(4)}}{k^2}
\end{eqnarray}
and
\begin{eqnarray}
\sum^{\infty}_{k=1} \left(\frac{1}{k^4}\sum^{k}_{i=1}\frac{h_i}{i^2}\right) = \zeta(4)\sum^{\infty}_{k=1}\frac{h_k}{k^2} + \sum^{\infty}_{k=1}\frac{h_k}{k^6} -
\sum^{\infty}_{k=1} \frac{h_k H_k^{(4)}}{k^2}
\end{eqnarray}
Eq.~(43) results as all Euler sums appearing in Eq.~(58) can be calculated in terms of zeta values.

As a second example we calculate the following expression \cite{bra21}:
\begin{eqnarray}
\sum^{\infty}_{k=1}\frac{1}{k^2}\left( \sum^{\infty}_{i=1}\frac{h_i^{(3)}}{i(k+i)} \right) = 2ln(2)\sum^{\infty}_{k=1}\frac{h_k^{(3)}}{k^3} + 
\frac{7}{4}\zeta(3)\sum^{\infty}_{k=1}\frac{h_k}{k^3} - \frac{3}{2}\zeta(2)\sum^{\infty}_{k=1}\frac{h_k^{(2)}}{k^3} + \sum^{\infty}_{k=1}\frac{1}{k^3}
\left( \sum^{k}_{i=1}\frac{H_{i-1}}{(2i-1)^3} \right)~.
\end{eqnarray}
Rearanging the summation order one gets:
\begin{eqnarray}
\sum^{\infty}_{k=1}\frac{1}{k^2}\left( \sum^{\infty}_{i=1}\frac{h_i^{(3)}}{i(k+i)} \right) = \zeta(2)\sum^{\infty}_{k=1}\frac{h_k^{(3)}}{k^2} -
\sum^{\infty}_{k=1} \frac{H_k h_k^{(3)}}{k^3}
\end{eqnarray}
With this it follows:
\begin{eqnarray}
\sum^{\infty}_{k=1} \frac{H_k h_k^{(3)}}{k^3} = \zeta(2)\sum^{\infty}_{k=1}\frac{h_k^{(3)}}{k^2} + \frac{3}{2}\zeta(2)\sum^{\infty}_{k=1}\frac{h_k^{(2)}}{k^3} -
2ln(2)\sum^{\infty}_{k=1}\frac{h_k^{(3)}}{k^3} - \frac{7}{4}\zeta(3)\sum^{\infty}_{k=1}\frac{h_k}{k^3} - \sum^{\infty}_{k=1}\frac{1}{k^3}
\left( \sum^{k}_{i=1}\frac{H_{i-1}}{(2i-1)^3} \right)~.
\end{eqnarray}
It remains to compute the last sum. By use of the following two valued help function:
\begin{eqnarray}
\sum^{\infty}_{i=1} \frac{h_i^{(2)}}{(i+k)^2} = \sum^{\infty}_{k=1}\frac{h_k^{(2)}}{k^2} - 4\zeta(2)h_k^{(2)} + 8ln(2)h_k^{(3)} + 4 
\sum^{k}_{i=1}\frac{H_{i-1}}{(2i-1)^3} + \sum^{k}_{i=1}\frac{H_{i-1}^{(2)}}{(2i-1)^2}~.
\end{eqnarray}
Now we calculate
\begin{eqnarray}
\sum^{\infty}_{k=1} \frac{1}{k^3} \left( \sum^{\infty}_{i=1} \frac{h_i^{(2)}}{(i+k)^2} \right) &=& \zeta(3)\sum^{\infty}_{k=1}\frac{h_k^{(2)}}{k^2} - 
4\zeta(2)\sum^{\infty}_{k=1}\frac{h_k^{(2)}}{k^3} + 8ln(2)\sum^{\infty}_{k=1}\frac{h_k^{(3)}}{k^3} + 4\sum^{\infty}_{k=1}\frac{1}{k^3}
\left( \sum^{k}_{i=1}\frac{H_{i-1}}{(2i-1)^3} \right) \nonumber \\ &+&  \sum^{\infty}_{k=1} \frac{H_k^{(2)} h_k^{(2)}}{k^3} -
\sum^{\infty}_{k=1} \frac{1}{k^3} \left( \sum^{k}_{i=1} \frac{h_k^{(2)}}{k^2} \right)
\end{eqnarray}
Rearanging the summation order in the last sum on the right side and in the sum appearing on the left side of Eq.~(65) it follows:
\begin{eqnarray}
\sum^{\infty}_{k=1}\frac{1}{k^3} \left( \sum^{k}_{i=1}\frac{H_{i-1}}{(2i-1)^3} \right) &=& \zeta(2)\sum^{\infty}_{k=1}\frac{h_k^{(2)}}{k^3} - 
2ln(2)\sum^{\infty}_{k=1}\frac{h_k^{(3)}}{k^3} + \frac{1}{4}\sum^{\infty}_{k=1}\frac{h_k^{(2)}}{k^5} - \frac{1}{4} \sum^{\infty}_{k=1} \frac{H_k^{(2)} h_k^{(2)}}{k^3} -
\frac{1}{4} \sum^{\infty}_{k=1} \frac{H_k^{(3)} h_k^{(2)}}{k^2} \nonumber \\ &+& \frac{1}{4}  \sum^{\infty}_{i=1} h_i^{(2)} \left( \sum^{\infty}_{k=1}
\frac{1}{k^3(i+k)^2} \right)~.
\end{eqnarray}
Explicit calculation of the last term in Eq.~(66) results in:
\begin{eqnarray}
\sum^{\infty}_{k=1}\frac{1}{k^3} \left( \sum^{k}_{i=1}\frac{H_{i-1}}{(2i-1)^3} \right) &=&  \frac{1}{4}\zeta(2)\sum^{\infty}_{k=1}\frac{h_k^{(2)}}{k^3} -
2ln(2)\sum^{\infty}_{k=1}\frac{h_k^{(3)}}{k^3} + \frac{1}{4}\sum^{\infty}_{k=1}\frac{h_k^{(2)}}{k^5} + \frac{1}{4}\zeta(3)\sum^{\infty}_{k=1}\frac{h_k^{(2)}}{k^2} -
\frac{1}{4} \sum^{\infty}_{k=1} \frac{H_k^{(3)} h_k^{(2)}}{k^2} \nonumber \\ &+& \frac{3}{4} \sum^{\infty}_{k=1} \frac{H_k h_k^{(2)}}{k^4}~.
\end{eqnarray}
Inserting Eq.~(67) in Eq.~(63) gives the following identity:
\begin{eqnarray}
\sum^{\infty}_{k=1} \frac{H_k h_k^{(3)}}{k^3} &=& \frac{5}{4}\zeta(2)\sum^{\infty}_{k=1}\frac{h_k^{(2)}}{k^3} + \zeta(2)\sum^{\infty}_{k=1}\frac{h_k^{(3)}}{k^2} -
\frac{7}{4}\zeta(3)\sum^{\infty}_{k=1}\frac{h_k}{k^3} - \frac{1}{4}\sum^{\infty}_{k=1}\frac{h_k^{(2)}}{k^5} - \frac{1}{4}\zeta(3)\sum^{\infty}_{k=1}\frac{h_k^{(2)}}{k^2} +
\frac{1}{4} \sum^{\infty}_{k=1} \frac{H_k^{(3)} h_k^{(2)}}{k^2} \nonumber \\ &-& \frac{3}{4} \sum^{\infty}_{k=1} \frac{H_k h_k^{(2)}}{k^4}~.
\end{eqnarray}
This way we receive a second identity between the members of the first family. Finally we present a calculational scheme that results in a third identity. We start with
the following two valued help function:
\begin{eqnarray}
\sum^{\infty}_{i=1} \frac{h_i}{(i+k)^2}	= \frac{7}{4}\zeta(3) + \zeta(2)h_k -4ln(2)h_k^{(2)} - 2\sum^{k}_{i=1}\frac{H_{i-1}}{(2i-1)^2} - H_k^{(2)}h_k +
\sum^{k}_{i=1} \frac{h_i}{i^2}~. 
\end{eqnarray} 
With this it follows:
\begin{eqnarray}
\sum^{\infty}_{k=1}\frac{1}{k^4} \left( \sum^{\infty}_{i=1} \frac{h_i}{(i+k)^2} \right) &=& \frac{7}{4}\zeta(3)\zeta(4) + \zeta(2)\sum^{\infty}_{k=1}\frac{h_k}{k^4} -
4ln(2)\sum^{\infty}_{k=1}\frac{h_k^{(2)}}{k^4} - \sum^{\infty}_{k=1}\frac{1}{k^4} \left( \sum^{k}_{i=1}\frac{H_{i-1}}{(2i-1)^2} \right)
- \sum^{\infty}_{k=1} \frac{H_k^{(2)} h_k^{(2)}}{k^4} \nonumber \\ &+& \sum^{\infty}_{k=1}\frac{1}{k^4} \left( \sum^{k}_{i=1} \frac{h_i}{i^2} \right)~. 
\end{eqnarray}
Rearanging the summation order in the last sum on the right side and in the sum appearing on the left side of Eq.~(70) it follows:
\begin{eqnarray}
\sum^{\infty}_{k=1}\frac{1}{k^4} \left( \sum^{k}_{i=1}\frac{H_{i-1}}{(2i-1)^2} \right) &=& \frac{7}{8}\zeta(3)\zeta(4) - \frac{3}{2}\zeta(2)\sum^{\infty}_{k=1}\frac{h_k}{k^4} +
\frac{1}{2}\sum^{\infty}_{k=1}\frac{h_k}{k^6} + \zeta(3)\sum^{\infty}_{k=1}\frac{h_k}{k^3} - 2ln(2)\sum^{\infty}_{k=1}\frac{h_k^{(2)}}{k^4} + 
2\sum^{\infty}_{k=1} \frac{H_k h_k}{k^5} \nonumber \\ &-&  \frac{1}{2}\sum^{\infty}_{k=1} \frac{H_k^{(4)} h_k}{k^2}~.
\end{eqnarray}
Furthermore one may calculate the following expression:
\begin{eqnarray}
\sum^{\infty}_{k=1}\frac{1}{k^3}\left( \sum^{\infty}_{i=1}\frac{h_i^{(2)}}{i(k+i)} \right) = -2ln(2)\sum^{\infty}_{k=1}\frac{h_k^{(2)}}{k^4} + 
\frac{3}{2}\zeta(2)\sum^{\infty}_{k=1}\frac{h_k}{k^4} - \sum^{\infty}_{k=1}\frac{1}{k^4} \left( \sum^{k}_{i=1}\frac{H_{i-1}}{(2i-1)^2} \right)~.
\end{eqnarray}
Rearanging the summation order one gets:
\begin{eqnarray}
\sum^{\infty}_{k=1} \frac{H_k h_k^{(2)}}{k^4} = \frac{3}{2}\zeta(2)\sum^{\infty}_{k=1}\frac{h_k}{k^4} - \zeta(3)\sum^{\infty}_{k=1}\frac{h_k^{(2)}}{k^2} +
\zeta(2)\sum^{\infty}_{k=1}\frac{h_k^{(2)}}{k^3} - -2ln(2)\sum^{\infty}_{k=1}\frac{h_k^{(2)}}{k^4} - \sum^{\infty}_{k=1}\frac{1}{k^4}
\left( \sum^{k}_{i=1}\frac{H_{i-1}}{(2i-1)^2} \right)~.
\end{eqnarray}
Combining finally Eq.~(71) and Eq.~(73) it follows:
\begin{eqnarray}
\sum^{\infty}_{k=1} \frac{H_k h_k^{(2)}}{k^4} = \frac{127}{32}\zeta(7) + \frac{31}{4}\zeta(2)\zeta(5) - \frac{41}{4}\zeta(3)\zeta(4) - 2\sum^{\infty}_{k=1} \frac{H_k h_k}{k^5} -
\frac{1}{2}\sum^{\infty}_{k=1} \frac{H_k^{(2)} h_k}{k^4}~.
\end{eqnarray}
These example may illustrate the complete procedure which results finally in the set of Eqs.~(42-49). What remains is to find a direct relation  between 
$\sum^{\infty}_{k=1} \frac{H_k^{(2)} h_k}{k^4}$ and $\sum^{\infty}_{k=1} \frac{H_k h_k}{k^5}$ in terms of zeta values and linear Euler sums and probably
an explicit calculation of $\sum^{\infty}_{k=1} \frac{H_k h_k}{k^5}$ in terms of  zeta values, linear Euler sums and $Li_7(1/2)$.  

\subsection{Second family}
All members of the second family can be calculated by use of lemma 2a from \cite{bra21} and with the Euler sums belonging to the first family. As a first example we present:
\begin{eqnarray}
\sum^{\infty}_{k=1} \frac{H_k H_k}{(2k-1)^5} &=& \frac{635}{64}\zeta(7) - \frac{315}{32}ln(2)\zeta(6) + \frac{315}{32}\zeta(6) + \frac{31}{8}\left(ln(2)\right)^2\zeta(5)
-\frac{31}{32}\zeta(2)\zeta(5) - \frac{31}{4}ln(2)\zeta(5) \nonumber \\ &-& \frac{93}{8}\zeta(5) - \frac{315}{64}\zeta(3)\zeta(4) +  \frac{15}{2}ln(2)\zeta(4) -
\frac{15}{8}\zeta(4) + \frac{49}{16}ln(2)\zeta(3)^2 - \frac{49}{16}\zeta(3)^2 + \frac{21}{4}\zeta(2)\zeta(3) \nonumber \\ &-& 7ln(2)\zeta(3) + \frac{7}{2}\zeta(3) +
6ln(2)\zeta(2) - 11\zeta(2) + 20ln(2) - 4\left(ln(2)\right)^2~.
\end{eqnarray}
The next example is:
\begin{eqnarray}
\sum^{\infty}_{k=1} \frac{H_k H_k^{(2)}}{(2k-1)^4} &=& \frac{889}{128}\zeta(7) + \frac{195}{16}ln(2)\zeta(6) - \frac{195}{16}\zeta(6) - \frac{155}{32}\zeta(2)\zeta(5) + 
\frac{155}{4}\zeta(5) + \frac{77}{32}\zeta(3)\zeta(4) - \frac{15}{2}ln(2)\zeta(4) \nonumber \\ &-& \frac{45}{8}\zeta(4) - \frac{63}{8}ln(2)\zeta(3)^2 + \frac{63}{8}\zeta(3)^2
- \frac{35}{2}\zeta(2)\zeta(3) + 14ln(2)\zeta(3) + \frac{7}{2}\zeta(3)- 20ln(2)\zeta(2) \nonumber \\ &+& 36\zeta(2) - 80ln(2) + 16\left(ln(2)\right)^2 +
2\sum^{\infty}_{k=1}\frac{h_k}{k^3} - ln(2)\sum^{\infty}_{k=1} \frac{h_k}{k^5} + \sum^{\infty}_{k=1}\frac{h_k}{k^5} - \frac{9}{8}\sum^{\infty}_{k=1} \frac{H_k^{(2)} h_k}{k^4} -
\nonumber \\ &-& \frac{1}{2}\sum^{\infty}_{k=1} \frac{H_k h_k}{k^5}~.
\end{eqnarray}

Also we found:
\begin{eqnarray}
\sum^{\infty}_{k=1} \frac{H_k H_k^{(3)}}{(2k-1)^3} &=& -\frac{3937}{256}\zeta(7) - \frac{135}{8}ln(2)\zeta(6) + \frac{135}{8}\zeta(6) + \frac{31}{32}\zeta(2)\zeta(5) - 
31\zeta(5) + \frac{517}{64}\zeta(3)\zeta(4) + \frac{45}{4}\zeta(4) \nonumber \\ &+& \frac{77}{8}ln(2)\zeta(3)^2 - \frac{77}{8}\zeta(3)^2 + 16\zeta(2)\zeta(3)
- 12ln(2)\zeta(3) - 40\zeta(3) + 36ln(2)\zeta(2) - 48\zeta(2) + 160ln(2) \nonumber \\ &-& 48\left(ln(2)\right)^2 - 2\sum^{\infty}_{k=1}\frac{h_k}{k^3}
+ 3ln(2)\sum^{\infty}_{k=1} \frac{h_k}{k^5} - 3\sum^{\infty}_{k=1}\frac{h_k}{k^5} + \frac{13}{8}\sum^{\infty}_{k=1} \frac{H_k^{(2)} h_k}{k^4} +
\frac{3}{2}\sum^{\infty}_{k=1} \frac{H_k h_k}{k^5}~.
\end{eqnarray}

Furthermore it follows:
\begin{eqnarray}
\sum^{\infty}_{k=1} \frac{H_k H_k^{(4)}}{(2k-1)^2} &=& \frac{127}{32}\zeta(7) + \frac{69}{8}ln(2)\zeta(6) - \frac{69}{8}\zeta(6) + \frac{93}{8}\zeta(2)\zeta(5) + 
\frac{31}{2}\zeta(5) - \frac{231}{16}\zeta(3)\zeta(4) - 2ln(2)\zeta(4) \nonumber \\ &+& 5\zeta(4) - \frac{7}{2}ln(2)\zeta(3)^2 + \frac{7}{2}\zeta(3)^2 - 7\zeta(2)\zeta(3)
+ 52\zeta(3) - 24ln(2)\zeta(2) + 24\zeta(2) - 160ln(2) \nonumber \\ &+& 64\left(ln(2)\right)^2 - \zeta(3)\sum^{\infty}_{k=1}\frac{h_k}{k^3} 
- 4ln(2)\sum^{\infty}_{k=1} \frac{h_k}{k^5} + 4\sum^{\infty}_{k=1}\frac{h_k}{k^5} - \frac{1}{2}\sum^{\infty}_{k=1} \frac{H_k^{(2)} h_k}{k^4} -
2\sum^{\infty}_{k=1} \frac{H_k h_k}{k^5}~,
\end{eqnarray}
and
\begin{eqnarray}
\sum^{\infty}_{k=1} \frac{H_k^{(2)} H_k^{(3)}}{(2k-1)^2} &=& \frac{9525}{64}\zeta(7) - \frac{217}{4}\zeta(2)\zeta(5) + 62\zeta(5) - 39\zeta(3)\zeta(4) - \frac{11}{2}\zeta(4)
- 31\zeta(2)\zeta(3) - 8ln(2)\zeta(3) + 30\zeta(3) \nonumber \\ &-& 24ln(2)\zeta(2) + 56\zeta(2) - 160ln(2) + \zeta(3)\sum^{\infty}_{k=1}\frac{h_k}{k^3} 
+ 16\sum^{\infty}_{k=1}\frac{h_k}{k^3} - 7\sum^{\infty}_{k=1} \frac{H_k^{(2)} h_k}{k^4}~.
\end{eqnarray}
Finally it follows:
\begin{eqnarray}
\sum^{\infty}_{k=1} \frac{H_k^{(2)} H_k^{(2)}}{(2k-1)^3} &=& -\frac{9017}{64}\zeta(7) + \frac{217}{4}\zeta(2)\zeta(5) - 93\zeta(5) + \frac{511}{16}\zeta(3)\zeta(4) + \frac{25}{2}\zeta(4)
+ 49\zeta(2)\zeta(3) \nonumber \\ &-& 10\zeta(3) + 24ln(2)\zeta(2) - 72\zeta(2) + 160ln(2) -16\sum^{\infty}_{k=1}\frac{h_k}{k^3} +7\sum^{\infty}_{k=1} \frac{H_k^{(2)} h_k}{k^4}~.
\end{eqnarray}
With this we have calculated all members of the second family concerning the order seven.

\subsection{Third family}
Again all members of the third family can be calculated by use of lemma 3a and lemma 3b from \cite{bra21} and with the Euler sums belonging to the first family.
As a first member of this family we show the following Euler sum:
\begin{eqnarray}
\sum^{\infty}_{k=1} \frac{H_k H_k^{(5)}}{k(2k-1)} &=& -6\zeta(7) - \frac{25}{2}\zeta(2)\zeta(5) - 2\left(ln(2)\right)^2\zeta(5) +4ln(2)\zeta(5) -12\zeta(5) + \frac{59}{4}\zeta(3)\zeta(4)
- 10\zeta(4) \nonumber \\ &+& 4\zeta(2)\zeta(3) - 32\zeta(3) + 128ln(2) -64\left(ln(2)\right)^2 + 2\zeta(3)\sum^{\infty}_{k=1}\frac{h_k}{k^3} + 4ln(2)\sum^{\infty}_{k=1} \frac{h_k}{k^5}
\nonumber \\ &-& 4\sum^{\infty}_{k=1} \frac{h_k}{k^5} +2\sum^{\infty}_{k=1} \frac{H_k h_k}{k^5}~.
\end{eqnarray}

Furthermore we found:
\begin{eqnarray}
\sum^{\infty}_{k=1} \frac{H_k^{(2)} H_k^{(4)}}{k(2k-1)} &=& \frac{7}{4}\zeta(7) + \frac{7}{2}ln(2)\zeta(6) - \frac{7}{2}\zeta(6) - \frac{15}{2}\zeta(2)\zeta(5) - 
32\zeta(5) + 6\zeta(3)\zeta(4) + 8ln(2)\zeta(4) -14\zeta(4) \nonumber \\ &+& 12\zeta(2)\zeta(3) - 40\zeta(3) + 32ln(2)\zeta(2) - 32\zeta(2) +128ln(2)~.
\end{eqnarray}

Finally it results:
\begin{eqnarray}
\sum^{\infty}_{k=1} \frac{H_k^{(3)} H_k^{(3)}}{k(2k-1)} &=& -\frac{4991}{16}\zeta(7) + \frac{237}{2}\zeta(2)\zeta(5) - 40\zeta(5) + \frac{307}{4}\zeta(3)\zeta(4) + 12\zeta(4)
+ 2ln(2)\zeta(3)^2 \nonumber \\ &-& 2\zeta(3)^2 + 16\zeta(2)\zeta(3) + 32ln(2)\zeta(3) - 16\zeta(3) - 32\zeta(2) + 128ln(2) - 4\zeta(3)\sum^{\infty}_{k=1}\frac{h_k}{k^3} \nonumber \\ &-&
32\sum^{\infty}_{k=1}\frac{h_k}{k^3} + 14\sum^{\infty}_{k=1} \frac{H_k^{(2)} h_k}{k^4}~.
\end{eqnarray}
With this we have calculated all members of the third family concerning the order seven.

\subsection{Fourth family}
Again all members of the fourth family can be calculated by use of lemma 4a and lemma 4b from \cite{bra21} and with the Euler sums belonging to the first family.
It follows first:
\begin{eqnarray}
\sum^{\infty}_{k=1} \frac{h_k^2}{(2k-1)^5} &=& \frac{127}{256}\zeta(7) + \frac{63}{64}ln(2)\zeta(6) - \frac{3}{128}\zeta(2)\zeta(5) + \frac{31}{32}\left(ln(2)\right)^2\zeta(5) 
-\frac{65}{256}\zeta(3)\zeta(4) - \frac{7}{32}ln(2)\zeta(3)^2  \nonumber \\ &-& \frac{1}{16}ln(2)\sum^{\infty}_{k=1}\frac{h_k}{k^5} -
\frac{1}{32}\sum^{\infty}_{k=1} \frac{H_k h_k}{k^5}~,
\end{eqnarray}

and for the second sum it results:
\begin{eqnarray}
\sum^{\infty}_{k=1} \frac{h_k h_k^{(2)}}{(2k-1)^4} &=& \frac{635}{1024}\zeta(7) + \frac{27}{32}ln(2)\zeta(6) - \frac{19}{256}\zeta(2)\zeta(5) - \frac{1}{8}\zeta(3)\zeta(4) +
\frac{7}{64}ln(2)\zeta(3)^2 - \frac{1}{64}\sum^{\infty}_{k=1} \frac{H_k^{(2)} h_k}{k^4}~. 
\end{eqnarray}
This results will be explicitly calculated as a showcase. Using Eq.~(160) from \cite{bra21} it follows:
\begin{eqnarray}
\sum^{\infty}_{k=1} \frac{1}{k^3} \left( \sum^{\infty}_{i=1} \frac{h_k^{(2)}}{(2i-1)(2i+2k-1)} \right) = \frac{7}{16}\zeta(3)\zeta(4) + \frac{3}{16}\sum^{\infty}_{k=1}\frac{H_{k-1}}{k^4} -
\frac{1}{8} \sum^{\infty}_{k=1} \frac{1}{k^4} \left( \sum^{k-1}_{i=1} \frac{h_i}{i^2} \right)~. 
\end{eqnarray}
By partial fraction decomposition of the left side and by rearanging the summation order fro the last sum of the right side one gets:
\begin{eqnarray}
\sum^{\infty}_{k=1} \frac{h_k h_k^{(2)}}{(2k-1)^4} &=& ln(2)\sum^{\infty}_{k=1} \frac{h_k^{(2)}}{(2k-1)^4} + \frac{1}{4}\zeta(2)\sum^{\infty}_{k=1} \frac{h_k^{(2)}}{(2k-1)^3}
- \frac{1}{8}\zeta(2)\sum^{\infty}_{k=1} \frac{h_k^{(2)}}{(2k-1)^2} - \frac{3}{64}\zeta(2)\zeta(5) - \frac{1}{32}\zeta(3)\zeta(4)  \nonumber \\ &+& \frac{1}{64}
\sum^{\infty}_{k=1} \frac{H_k^{(4)} h_k}{k^2}~.
\end{eqnarray}
As all linear sums in Eq.~(87) are known in terms of zeta values Eq.~(85) results. 

Furthermore it follows:
\begin{eqnarray}
\sum^{\infty}_{k=1} \frac{h_k h_k^{(3)}}{(2k-1)^3} &=& \frac{5715}{4096}\zeta(7) + \frac{63}{128}ln(2)\zeta(6) - \frac{283}{512}\zeta(2)\zeta(5) - 
\frac{43}{1024}\zeta(3)\zeta(4) + \frac{49}{128}ln(2)\zeta(3)^2  \nonumber \\ &-& \frac{7}{128}\sum^{\infty}_{k=1} \frac{H_k^{(2)} h_k}{k^4}~.
\end{eqnarray}

and
\begin{eqnarray}
\sum^{\infty}_{k=1} \frac{h_k h_k^{(4)}}{(2k-1)^2} &=& \frac{351}{256}ln(2)\zeta(6) + \frac{29}{128}\zeta(2)\zeta(5) + \frac{49}{512}\zeta(3)\zeta(4) - \frac{7}{64}ln(2)\zeta(3)^2
+ \frac{1}{64}\sum^{\infty}_{k=1} \frac{H_k^{(2)} h_k}{k^4}~.
\end{eqnarray}

For the last two Euler sums of this family we found:
\begin{eqnarray}
\sum^{\infty}_{k=1} \frac{h_k^{(2)} h_k^{(3)}}{(2k-1)^2} &=& \frac{635}{2048}\zeta(7) - \frac{15}{256}\zeta(2)\zeta(5) + \frac{105}{128}\zeta(3)\zeta(4)~,
\end{eqnarray}

and
\begin{eqnarray}
\sum^{\infty}_{k=1} \frac{h_k^{(2)} h_k^{(2)}}{(2k-1)^3} &=& -\frac{127}{1024}\zeta(7) + \frac{15}{128}\zeta(2)\zeta(5) + \frac{195}{256}\zeta(3)\zeta(4)~.
\end{eqnarray}
With this we have calculated all members of the fourth family concerning the order seven.

\subsection{Fifth family}
The calculational scheme for the members of the fifth family is more tricky. As a first member we show the following Euler sum: 
\begin{eqnarray}
\sum^{\infty}_{k=1} \frac{H_k h_k}{(2k-1)^5} &=& \frac{1905}{2048}\zeta(7) + \frac{189}{128}ln(2)\zeta(6) + \frac{63}{64}\zeta(6) - \frac{31}{256}\zeta(2)\zeta(5) - 
\frac{31}{16}\left(ln(2)\right)^2\zeta(5) + \frac{31}{16}ln(2)\zeta(5)  - \frac{31}{32}\zeta(5) \nonumber \\ &-&
\frac{167}{512}\zeta(3)\zeta(4) - \frac{15}{8}ln(2)\zeta(4) + \frac{15}{16}\zeta(4) - \frac{35}{64}ln(2)\zeta(3)^2 - \frac{7}{32}\zeta(3)^2 + \frac{3}{16}\zeta(2)\zeta(3) +
\frac{7}{4}ln(2)\zeta(3) \nonumber \\ &-& \frac{7}{8}\zeta(3) - \frac{3}{2}ln(2)\zeta(2) + \zeta(2) - \frac{1}{4}\sum^{\infty}_{k=1}\frac{h_k}{k^3} +
\frac{1}{16}ln(2)\sum^{\infty}_{k=1}\frac{h_k}{k^5} - \frac{1}{16}\sum^{\infty}_{k=1}\frac{h_k}{k^5} + \frac{7}{64}\sum^{\infty}_{k=1} \frac{H_k^{(2)} h_k}{k^4}~.
\end{eqnarray}
This sum results from an explicit calculation of the following expression, namely by performing a corresponding partial fraction decomposition :
\begin{eqnarray}
\sum^{\infty}_{k=1}\frac{H_k}{2k-1} \left( \sum^{\infty}_{i=1} \frac{1}{i^4 (2i+2k-1)} \right)~.
\end{eqnarray}
The three Euler sums presented within Eq.~(102), Eq.~(104) and Eq.~(107) can be calculated explicitly by use of Eq.~(246) from \cite{ade16,bra21} and by use of members of the
first family. For the three Euler sums presented within Eq.~(99), Eq.~(103) and Eq.~(109) another approach is needed. For example it stands:
\begin{eqnarray}
\sum^{\infty}_{k=1} \frac{H_k h_k^{(2)}}{(2k-1)^4} = \frac{3}{4}\zeta(2) \sum^{\infty}_{k=1} \frac{H_k}{(2k-1)^4} - \sum^{\infty}_{k=1} \frac{H_k}{(2k-1)^4} \left( 
\sum^{\infty}_{i=1} \frac{1}{(2i+2k-1)^2} \right)~.
\end{eqnarray}
By rearanging the summation order and performing a partial fraction decomposition it follows:
\begin{eqnarray}
\sum^{\infty}_{k=1} \frac{H_k h_k^{(2)}}{(2k-1)^4} &=& \frac{3}{4}\zeta(2) \sum^{\infty}_{k=1} \frac{H_k}{(2k-1)^4} - \sum^{\infty}_{i=1} \Big( \sum^{\infty}_{k=1} \Big( 
\frac{1}{4i^2}\frac{H_k}{(2k-1)^4} - \frac{1}{4i^3}\frac{H_k}{(2k-1)^3} + \frac{3}{16i^2}\frac{H_k}{(2k-1)^2} \nonumber \\ &-& \frac{1}{4i^4}\frac{H_k}{(2k-1)(2i+2k-1)} +
\frac{1}{16i^2}\frac{H_k}{(2i+2k-1)^2} \Big ) \Big )~.
\end{eqnarray}
Analogously to lemma 4b from \cite{bra21} it stands:
\begin{eqnarray}
\sum^{\infty}_{k=1} \frac{H_k}{(2i+2k-1)^2} = \frac{H_i}{(2i+1)^2} + \frac{1}{4}\zeta(2)h_{i-1} + 2ln(2)h_{i-1}^{(2)} - 2h_{i-1}^{(3)} -2h_{i-1} h_{i-1}^{(2)}~. 
\end{eqnarray}
Furthermore, analogously to lemma 2a from \cite{bra21} it stands: 
\begin{eqnarray}
\sum^{\infty}_{k=1} \frac{H_k}{(2k-1)(2i+2k-1)} = \frac{2ln(2)}{2i-1} - ln(2)\frac{h_i}{i} - \frac{h_i}{i(2i-1)} + \frac{h_i^{(2)}}{2i} + \frac{h_i^{2}}{2i}~.
\end{eqnarray}
With this two valued help functions we get:
\begin{eqnarray}
\sum^{\infty}_{k=1} \frac{H_k h_k^{(2)}}{(2k-1)^4} &=& \frac{1}{2}\zeta(2)\sum^{\infty}_{k=1} \frac{H_k}{(2k-1)^4} - \frac{1}{4}\zeta(3)\sum^{\infty}_{k=1} \frac{H_k}{(2k-1)^3}
+ \frac{1}{4}\zeta(4)\sum^{\infty}_{k=1}\frac{H_k}{(2k-1)^2} - \frac{1}{16}\zeta(4)\sum^{\infty}_{k=1} \frac{1}{k(2k-1)^2} \nonumber \\ &-& \frac{1}{2}ln(2)\sum^{\infty}_{k=1}
\frac{1}{k^4(2k-1)} + \frac{1}{4}ln(2)\sum^{\infty}_{k=1} \frac{h_k}{k^5} + \frac{1}{4}\sum^{\infty}_{k=1} \frac{h_k}{k^5(2k-1)} - \frac{1}{8}\sum^{\infty}_{k=1} \frac{h_k^{(2)}}{k^5}
- \frac{1}{8}\sum^{\infty}_{k=1} \frac{h_k^2}{k^5} \nonumber \\ &+& \frac{3}{32}\zeta(2)\sum^{\infty}_{k=1} \frac{h_{k-1}}{k^4} + \frac{1}{8}ln(2)\sum^{\infty}_{k=1} \frac{h_{k-1}^{(2)}}{k^4}
- \frac{1}{8}\sum^{\infty}_{k=1} \frac{h_{k-1}^{(3)}}{k^4} - \frac{1}{8}\sum^{\infty}_{k=1} \frac{h_{k-1} h_{k-1}^{(2)} }{k^4}~.
\end{eqnarray}
This way the corresponding Euler sum, belonging to the fifth family, depends on two Euler sums from family eight. Explicit calculation results in:
\begin{eqnarray}
\sum^{\infty}_{k=1} \frac{H_k h_k^{(2)}}{(2k-1)^4} &=& \frac{4953}{1024}\zeta(7) - \frac{27}{16}ln(2)\zeta(6) + \frac{27}{16}\zeta(6) -  \frac{93}{128}\zeta(2)\zeta(5)
- \frac{31}{32}\zeta(5) - \frac{435}{256}\zeta(3)\zeta(4) + \frac{75}{32}\zeta(4) \nonumber \\ &-& \frac{7}{32}ln(2)\zeta(3)^2 + \frac{7}{32}\zeta(3)^2 - \frac{9}{16}\zeta(2)\zeta(3) -
\frac{21}{8}\zeta(3) + \frac{3}{2}ln(2)\zeta(2)~.
\end{eqnarray}

Furthermore we found:
\begin{eqnarray}
\sum^{\infty}_{k=1} \frac{H_k^{(2)} h_k}{(2k-1)^4} &=& -\frac{8255}{512}\zeta(7) - \frac{195}{32}ln(2)\zeta(6) + \frac{403}{64}\zeta(2)\zeta(5) + \frac{31}{16}\zeta(5) +
\frac{519}{128}\zeta(3)\zeta(4) + \frac{15}{4}ln(2)\zeta(4) - \frac{15}{4}\zeta(4) \nonumber \\ &+& \frac{63}{16}ln(2)\zeta(3)^2 - \frac{3}{8}\zeta(2)\zeta(3) -
7ln(2)\zeta(3) + 7\zeta(3) + 9ln(2)\zeta(2) - 8\zeta(2) + \sum^{\infty}_{k=1}\frac{h_k}{k^3} +
\frac{1}{2}ln(2)\sum^{\infty}_{k=1}\frac{h_k}{k^5} \nonumber \\ &-& \frac{1}{4}\sum^{\infty}_{k=1} \frac{H_k^{(2)} h_k}{k^4} + \frac{1}{2}\sum^{\infty}_{k=1} \frac{H_k h_k}{k^5}~.
\end{eqnarray}

The fourth sum results to:
\begin{eqnarray}
\sum^{\infty}_{k=1} \frac{H_k^{(3)} h_k}{(2k-1)^3} &=& \frac{14605}{256}\zeta(7) + \frac{135}{16}ln(2)\zeta(6) - \frac{837}{32}\zeta(2)\zeta(5) - 
- \frac{551}{64}\zeta(3)\zeta(4) + \frac{15}{4}\zeta(4) - \frac{77}{16}ln(2)\zeta(3)^2 + 7ln(2)\zeta(3) \nonumber \\ &-& 21\zeta(3) -18ln(2)\zeta(2) + 24\zeta(2)
- \frac{1}{8}\zeta(3)\sum^{\infty}_{k=1}\frac{h_k}{k^3} - 2\sum^{\infty}_{k=1} \frac{h_k}{k^3} - \frac{3}{2}ln(2)\sum^{\infty}_{k=1}\frac{h_k}{k^5}
+ \frac{1}{8}\sum^{\infty}_{k=1} \frac{H_k^{(2)} h_k}{k^4} \nonumber \\ &-& \frac{3}{2}\sum^{\infty}_{k=1} \frac{H_k h_k}{k^5}~.
\end{eqnarray}

It follows also:
\begin{eqnarray}
\sum^{\infty}_{k=1} \frac{H_k^{(2)} h_k^{(2)}}{(2k-1)^3} &=& -\frac{4699}{256}\zeta(7) + \frac{279}{64}\zeta(2)\zeta(5) + \frac{31}{16}\zeta(5) +
\frac{281}{32}\zeta(3)\zeta(4) - \frac{195}{16}\zeta(4) + \frac{9}{8}\zeta(2)\zeta(3) + \frac{63}{4}\zeta(3) \nonumber \\ &-& 9ln(2)\zeta(2)
+ \sum^{\infty}_{k=1}\frac{h_k}{k^3} - \frac{1}{4}\sum^{\infty}_{k=1} \frac{H_k^{(2)} h_k}{k^4}~,
\end{eqnarray}

and we found:
\begin{eqnarray}
\sum^{\infty}_{k=1} \frac{H_k h_k^{(3)}}{(2k-1)^3} &=& \frac{8509}{1024}\zeta(7) - \frac{63}{64}ln(2)\zeta(6) + \frac{63}{64}\zeta(6) - \frac{155}{64}\zeta(2)\zeta(5) - 
\frac{31}{32}\zeta(5) - \frac{553}{256}\zeta(3)\zeta(4) + \frac{75}{32}\zeta(4) \nonumber \\ &-& \frac{49}{64}ln(2)\zeta(3)^2 + \frac{49}{64}\zeta(3)^2 -
\frac{3}{4}\zeta(2)\zeta(3) - \frac{21}{12}ln(2)\zeta(3) + \frac{7}{32}\zeta(3)\sum^{\infty}_{k=1}\frac{h_k}{k^3} + \frac{1}{4}\sum^{\infty}_{k=1} \frac{h_k}{k^3}
\nonumber \\ &-& \frac{7}{32}\sum^{\infty}_{k=1} \frac{H_k^{(2)} h_k}{k^4}~.
\end{eqnarray}

Furthermore it results:
\begin{eqnarray}
\sum^{\infty}_{k=1} \frac{H_k^{(2)} h_k^{(3)}}{(2k-1)^2} &=& -\frac{22479}{512}\zeta(7) + \frac{961}{64}\zeta(2)\zeta(5) + \frac{217}{16}\zeta(5) + \frac{117}{8}\zeta(3)\zeta(4)
- \frac{75}{8}\zeta(4) \nonumber \\ &-& \frac{15}{4}\zeta(2)\zeta(3) + 7ln(2)\zeta(3) - \frac{7}{8}\zeta(3)\sum^{\infty}_{k=1}\frac{h_k}{k^3} - \sum^{\infty}_{k=1} \frac{h_k}{k^3}
+ \frac{9}{8}\sum^{\infty}_{k=1} \frac{H_k^{(2)} h_k}{k^4}~.
\end{eqnarray}

For the three Euler sums presented by Eq.~(100), Eq.~(101) and Eq.~(106) a third type of approach has been found to be successful. This will be shown in the following. We start
with the identity:
\begin{eqnarray}
\sum^{\infty}_{i=1} \frac{1}{i^2} \left( \sum^{\infty}_{k=1} \frac{h_{k+i}}{(k+i)^4} \right) &=& \sum^{\infty}_{k=1} \frac{1}{k^4} \left( \sum^{\infty}_{i=1} \frac{h_{k+i}}{i^2} \right) -
4\sum^{\infty}_{k=1} \frac{1}{k^4} \left( \sum^{\infty}_{i=1} \frac{h_{k+i}}{i(k+i)} \right) +3 \sum^{\infty}_{k=1} \frac{1}{k^4} \left( \sum^{\infty}_{i=1} \frac{h_{k+i}}{(k+i)^2} \right)
\nonumber \\ &+& 2\sum^{\infty}_{k=1} \frac{1}{k^3} \left( \sum^{\infty}_{i=1} \frac{h_{k+i}}{(i+k)^3} \right)~.
\end{eqnarray}
Analogously to Eq.~(242) in \cite{bra21} it stands:
\begin{eqnarray}
\sum^{\infty}_{i=1} \frac{h_{k+i}}{i(k+i)} = -2ln(2)\frac{h_k}{k} + \frac{h_k^{(2)}}{k} + \frac{h_k^2}{k} + \frac{1}{k}\sum^{k}_{i=1} \frac{h_i}{i}
\end{eqnarray}
and furthermore it stands:
\begin{eqnarray}
\frac{h_{k+i}}{i^2} = \frac{7}{4}\zeta(3) + \zeta(2)h_k + 4ln(2)h_k^{(2)} - 4\sum^{k}_{i=1} \frac{h_i}{(2i-1)^2}~.
\end{eqnarray}
Using these two identities we finally arrive at the following expression for Eq.~(109)
\begin{eqnarray}
&& 8\sum^{\infty}_{k=1} \frac{H_k h_k}{k^5} - 2\sum^{\infty}_{k=1} \frac{H_k^{(3)} h_k}{k^3} -3\sum^{\infty}_{k=1} \frac{H_k^{(4)} h_k}{k^2} + 4\sum^{\infty}_{k=1} \frac{h_k^2}{k^5} =
\frac{7}{4}\zeta(3)\zeta(4) + 4\zeta(4)\sum^{\infty}_{k=1}\frac{h_k}{k^2} - 4\zeta(3)\sum^{\infty}_{k=1}\frac{h_k}{k^3} + 5\zeta(2)\sum^{\infty}_{k=1}\frac{h_k}{k^4} \nonumber \\ &-& 
4ln(2)\sum^{\infty}_{k=1}\frac{h_k^{(2)}}{k^4} - 4\sum^{\infty}_{k=1}\frac{h_k^{(2)}}{k^5} - 5\sum^{\infty}_{k=1}\frac{h_k}{k^6} - 4\zeta(4)\sum^{\infty}_{k=1}\frac{h_k}{(2k-1)^2} -
4\sum^{\infty}_{k=1}\frac{h_k}{k^4(2k-1)^2} + 4\sum^{\infty}_{k=1} \frac{H_k^{(4)} h_k}{(2k-1)^2}~.
\end{eqnarray}
From the last expression Eq.~(109) follows after some algebra, where members of the first and eight family are involved.
\begin{eqnarray}
\sum^{\infty}_{k=1} \frac{H_k^{(4)} h_k}{(2k-1)^2} &=& -\frac{8255}{182}\zeta(7) - \frac{69}{16}ln(2)\zeta(6) + \frac{279}{8}\zeta(2)\zeta(5) + 
\frac{31}{4}\zeta(5) + \frac{49}{16}\zeta(3)\zeta(4) + \frac{7}{4}ln(2)\zeta(3)^2 - \frac{7}{2}\zeta(2)\zeta(3) \nonumber \\ &+& 28\zeta(3) + 12ln(2)\zeta(2) - 32\zeta(2) +
\frac{1}{2}\zeta(3)\sum^{\infty}_{k=1}\frac{h_k}{k^3} + 4\sum^{\infty}_{k=1} \frac{h_k}{k^3} + 2ln(2)\sum^{\infty}_{k=1}\frac{h_k}{k^5} - 
\sum^{\infty}_{k=1} \frac{H_k^{(2)} h_k}{k^4} \nonumber \\ &+& 2\sum^{\infty}_{k=1} \frac{H_k h_k}{k^5}~.
\end{eqnarray}

Next to last it follows:
\begin{eqnarray}
\sum^{\infty}_{k=1} \frac{H_k^{(3)} h_k^{(2)}}{(2k-1)^2} &=& \frac{4953}{256}\zeta(7) - \frac{31}{4}\zeta(2)\zeta(5) - \frac{31}{2}\zeta(5) - \frac{197}{32}\zeta(3)\zeta(4) +
\frac{165}{8}\zeta(4) + \frac{35}{4}\zeta(2)\zeta(3) \nonumber \\ &-& \frac{63}{2}\zeta(3) + 18ln(2)\zeta(2) - 4\sum^{\infty}_{k=1} \frac{h_k}{k^3} + 
\frac{7}{4}\sum^{\infty}_{k=1} \frac{H_k^{(2)} h_k}{k^4}~.
\end{eqnarray}

Finally we present:
\begin{eqnarray}
\sum^{\infty}_{k=1} \frac{H_k h_k^{(4)}}{(2k-1)^2} &=& \frac{4953}{1024}\zeta(7) - \frac{351}{128}ln(2)\zeta(6) + \frac{351}{128}\zeta(6) - \frac{93}{128}\zeta(2)\zeta(5) - 
\frac{155}{32}\zeta(5) - \frac{75}{64}\zeta(3)\zeta(4) + \frac{15}{8}ln(2)\zeta(4) \nonumber \\ &+& \frac{7}{32}ln(2)\zeta(3)^2 - \frac{7}{32}\zeta(3)^2 + \frac{9}{8}\zeta(2)\zeta(3)~.
\end{eqnarray}
With this we have calculated all members of the fifth family to the order 7.

\subsection{Sixth family}
As for the fifth family the situation is similar for the sixth and seventh family. The explicit calculation of different members needs different approaches. For example, for the first
Euler sum of order seven belonging to the sixth family it follows:
\begin{eqnarray}
\sum^{\infty}_{k=1} \frac{h_k h_k^{(5)}}{k(2k-1)} &=& -\frac{4191}{1024}\zeta(7) + \frac{63}{64}ln(2)\zeta(6) + \frac{73}{32}\zeta(2)\zeta(5) - \frac{31}{32}\left(ln(2)\right)^2\zeta(5)
+ \frac{9}{8}\zeta(3)\zeta(4) + \frac{7}{32}ln(2)\zeta(3)^2 \nonumber \\ &-& \frac{7}{32}\zeta(3)\sum^{\infty}_{k=1}\frac{h_k}{k^3} + \frac{1}{16}ln(2)\sum^{\infty}_{k=1}\frac{h_k}{k^5} 
+ \frac{7}{32}\sum^{\infty}_{k=1} \frac{H_k^{(2)} h_k}{k^4} +  \frac{1}{32}\sum^{\infty}_{k=1} \frac{H_k h_k}{k^5}~.
\end{eqnarray}
This result has been achieved by using the following two valued help function:
\begin{eqnarray}
\sum^{\infty}_{k=1} \frac{h_{k+i}}{k(2k-1)} = \zeta(2) + 2ln(2)h_i -  \sum^{i}_{k=1} \frac{h_{k}}{k(2k-1)}~.
\end{eqnarray}
With this it follows:
\begin{eqnarray}
\sum^{\infty}_{i=1} \frac{1}{(2i-1)^5} \left( \sum^{\infty}_{k=1} \frac{h_{k+i}}{k(2k-1)} \right) = \frac{31}{32}\zeta(2)\zeta(5) + 2ln(2)\sum^{\infty}_{i=1}\frac{h_i}{(2i-1)^5} -
\sum^{\infty}_{i=1} \frac{1}{(2i-1)^5} \left( \sum^{i}_{k=1} \frac{h_{k}}{k(2k-1)} \right) ~.
\end{eqnarray}
Rearanging the summation order on left side and using then the two valued help function:
\begin{eqnarray}
\frac{h_{k+i}}{(2k-1)^5} =  \frac{h_k}{(2k-1)^5} + \frac{31}{64}\zeta(5)H_i - \frac{15}{64}\zeta(5)H_i^{(2)} + \frac{7}{64}\zeta(5)H_i^{(3)} - \frac{3}{64}\zeta(5)H_i^{(4)} +
\frac{1}{32}\sum^{i}_{k=1} \frac{h_k}{k^5} 
\end{eqnarray}
it follows:
\begin{eqnarray}
\sum^{\infty}_{k=1} \frac{h_k h_k^{(5)}}{k(2k-1)} &=& \sum^{\infty}_{k=1} \frac{h_k}{k(2k-1)^6} + \frac{31}{64}\zeta(5)\sum^{\infty}_{k=1}\frac{H_k}{k(2k-1)} - 
\frac{15}{64}\zeta(4)\sum^{\infty}_{k=1}\frac{H_k^{(2)}}{k(2k-1)} + \frac{7}{64}\zeta(3)\sum^{\infty}_{k=1}\frac{H_k^{(3)}}{k(2k-1)} \nonumber \\ &-& \frac{3}{64}\zeta(2)\sum^{\infty}_{k=1}
\frac{H_k^{(4)}}{k(2k-1)} + \frac{1}{16}ln(2)\sum^{\infty}_{k=1} \frac{h_k}{k^5} + \frac{1}{32}\sum^{\infty}_{k=1}\frac{h_k}{k^6 (2k-1)} - \frac{1}{16}\sum^{\infty}_{k=1} \frac{h_k^2}{k^5}
+ \frac{1}{32}\sum^{\infty}_{k=1} \frac{H_k h_k}{k^5}~. 
\end{eqnarray}
As one can observe, members of the first and of the eighth family are needed for an explicit calculation. Furthermore, for the second sum we found:
\begin{eqnarray}
\sum^{\infty}_{k=1} \frac{h_k^{(3)} h_k^{(3)}}{k(2k-1)} &=& \frac{6223}{1024}\zeta(7) - \frac{277}{128}\zeta(2)\zeta(5) + \frac{17}{256}\zeta(3)\zeta(4)
- \frac{49}{32}ln(2)\zeta(3)^2 + \frac{7}{16}\zeta(3)\sum^{\infty}_{k=1}\frac{h_k}{k^3} - \frac{7}{32}\sum^{\infty}_{k=1} \frac{H_k^{(2)} h_k}{k^4}~.
\end{eqnarray}
In this case again \cite{ade16} is very helpful. We start the calculation with:
\begin{eqnarray}
\sum^{k}_{i=1} \frac{h_i^{(3)}}{(2i-1)} = h_k h_k^{(3)} + h_k^{(4)} - \sum^{k}_{i=1} \frac{h_i}{(2i-1)^3}~.
\end{eqnarray}
From this it follows:
\begin{eqnarray}
\sum^{\infty}_{k=1} \frac{h_k^{(3)} h_k^{(3)}}{k(2k-1)} &=& \frac{7}{8}\zeta(3)\sum^{\infty}_{k=1} \frac{h_k^{(3)}}{k(2k-1)} + \sum^{\infty}_{k=1} \frac{h_k^{(3)}}{k(2k-1)^4} -
2\sum^{\infty}_{k=1} \frac{h_k h_k^{(3)}}{(2k-1)^3} - 2\sum^{\infty}_{k=1} \frac{h_k^{(4)}}{(2k-1)^3} + \sum^{\infty}_{k=1} \frac{H_k h_k^{(3)}}{(2k-1)^3} \nonumber \\ &+&
2\sum^{\infty}_{k=1} \frac{1}{(2k-1)^3} \left( \sum^{k}_{i=1} \frac{h_i}{(2i-1)^3} \right) - \sum^{\infty}_{k=1} \frac{1}{(2k-1)^3} \left( \sum^{k}_{i=1} \frac{H_{i-1}}{(2i-1)^3} \right)
\end{eqnarray}
Rearanging the summation order in the last two sums on the right side we get:
\begin{eqnarray}
\sum^{\infty}_{k=1} \frac{h_k^{(3)} h_k^{(3)}}{k(2k-1)} &=& \frac{7}{8}\zeta(3)\sum^{\infty}_{k=1} \frac{h_k^{(3)}}{k(2k-1)} + \sum^{\infty}_{k=1} \frac{h_k^{(3)}}{k(2k-1)^4} -
2\sum^{\infty}_{k=1} \frac{h_k^{(4)}}{(2k-1)^3} + \frac{7}{4}\zeta(3)\sum^{\infty}_{k=1} \frac{h_k}{(2k-1)^3} + 2\sum^{\infty}_{k=1} \frac{h_k}{(2k-1)^6} \nonumber \\ &-&
\frac{7}{8}\zeta(3)\sum^{\infty}_{k=1} \frac{H_{k-1}}{(2k-1)^3} - \sum^{\infty}_{k=1} \frac{H_{k-1}}{(2k-1)^6} - \sum^{\infty}_{k=1}\frac{h_k^{(3)}}{k(2k-1)^3}
-4\sum^{\infty}_{k=1} \frac{h_k h_k^{(3)}}{(2k-1)^3} + 2\sum^{\infty}_{k=1} \frac{H_k h_k^{(3)}}{(2k-1)^3}~.
\end{eqnarray}
As one can observe, members of the fourth and of the fifth family are needed for an explicit calculation. Finally it results:
\begin{eqnarray}
\sum^{\infty}_{k=1} \frac{h_k^{(2)} h_k^{(4)}}{k(2k-1)} &=& \frac{889}{256}\zeta(7) - \frac{315}{128}ln(2)\zeta(6) + \frac{45}{128}\zeta(2)\zeta(5) - \frac{45}{64}\zeta(3)\zeta(4)~.
\end{eqnarray}
In this case again \cite{ade16} is very helpful. We start the calculation with:
\begin{eqnarray}
\sum^{k}_{i=1} \frac{h_i^{(2)}}{(2i-1)^4} = h_k^{(2)} h_k^{(4)} + h_k^{(6)} - \sum^{k}_{i=1} \frac{h_i^{(4)}}{(2i-1)^2}~. 
\end{eqnarray}
From this it follows:
\begin{eqnarray}
\sum^{\infty}_{k=1} \frac{1}{k(2k-1)} \left( \sum^{k}_{i=1} \frac{h_i^{(2)}}{(2i-1)^4} \right) &=& \sum^{\infty}_{k=1} \frac{h_k^{(2)} h_k^{(4)}}{k(2k-1)} + \sum^{\infty}_{k=1}
\frac{h_k^{(6)}}{k(2k-1)} - \sum^{\infty}_{k=1} \frac{1}{k(2k-1)} \left( \sum^{k}_{i=1} \frac{h_i^{(4)}}{(2i-1)^2} \right)~.
\end{eqnarray}
The algebra needed for an explicit calculation of Eq.~(121) is analogous to the computation shown for Eq.~(117). Again, one can observe, members of the fifth family are needed for
an explicit calculation.  With this we have calculated all members of the sixth family to the order seven.

\subsection{Seventh family}
As a first member we show the following Euler sum:
\begin{eqnarray}
\sum^{\infty}_{k=1} \frac{H_k h_k^{(5)}}{k(2k-1)} &=& -\frac{8255}{1024}\zeta(7) - \frac{189}{16}ln(2)\zeta(6) + \frac{441}{64}\zeta(6) + \frac{775}{128}\zeta(2)\zeta(5) + 
\frac{31}{8}\left(ln(2)\right)^2\zeta(5) - \frac{31}{8}ln(2)\zeta(5) \nonumber \\ &+&  \frac{223}{256}\zeta(3)\zeta(4) + \frac{21}{8}ln(2)\zeta(3)^2 - \frac{35}{32}\zeta(3)^2
- \frac{7}{16}\zeta(3)\sum^{\infty}_{k=1}\frac{h_k}{k^3} - \frac{1}{8}ln(2)\sum^{\infty}_{k=1}\frac{h_k}{k^5} + \frac{1}{8}\sum^{\infty}_{k=1}\frac{h_k}{k^5}  \nonumber \\ &+&
\frac{7}{32}\sum^{\infty}_{k=1} \frac{H_k^{(2)} h_k}{k^4}~.
\end{eqnarray}
The calculational procedure is as follows. We start with the identity \cite{ade16}:
\begin{eqnarray}
&&\sum^{k}_{i=1} \frac{h_i}{(2i-1)^5} = h_k h_k^{(5)} + h_k^{(6)} - \sum^{k}_{i=1} \frac{h_i^{(5)}}{(2i-1)} = \nonumber \\ &&
h_k h_k^{(5)} + h_k^{(6)} - \frac{1}{2}\sum^{k}_{i=1} \frac{h_i^{(5)}}{i(2i-1)} - \frac{1}{2}H_k h_k^{(5)} + \frac{1}{2}\sum^{k}_{i=1} \frac{H_i}{(2i-1)^5} -
\frac{1}{2}\sum^{k}_{i=1} \frac{1}{i(2i-1)^5}~.
\end{eqnarray}
From this we get:
\begin{eqnarray}
\sum^{\infty}_{k=1} \frac{1}{k(2k-1)} \left( \sum^{k}_{i=1} \frac{h_i}{(2i-1)^5} \right)&=& \sum^{\infty}_{k=1} \frac{h_k h_k^{(5)}}{k(2k-1)}  + \sum^{\infty}_{k=1} \frac{h_k^{(6)}}{k(2k-1)}
- \frac{1}{2}\sum^{\infty}_{k=1} \frac{H_k h_k^{(5)}}{k(2k-1)} + \frac{1}{2} \sum^{\infty}_{k=1} \frac{1}{k(2k-1)} \left( \sum^{k}_{i=1} \frac{H_i}{(2i-1)^5} \right) \nonumber \\ &-&
\frac{1}{2} \sum^{\infty}_{k=1} \frac{1}{k(2k-1)} \left( \sum^{k}_{i=1} \frac{h_i^{(5)}}{i(2i-1)} \right) - \frac{1}{2} \sum^{\infty}_{k=1} \frac{1}{k(2k-1)} \left( \sum^{k}_{i=1}
\frac{1}{i(2i-1)^5} \right)~.
\end{eqnarray}
By rearanging the summation order in the last two sums of Eq.~(126) it follows:
\begin{eqnarray}
\sum^{\infty}_{k=1} \frac{H_k h_k^{(5)}}{k(2k-1)} &=& \sum^{\infty}_{k=1} \frac{h_k}{k(2k-1)^5} - \frac{1}{2}\sum^{\infty}_{k=1} \frac{H_k}{k(2k-1)^5} - 2ln(2)
\sum^{\infty}_{k=1} \frac{h_k}{(2k-1)^5} - \sum^{\infty}_{k=1} \frac{h_k}{k(2k-1)^6} + 2\sum^{\infty}_{k=1} \frac{h_k^2}{(2k-1)^5} \nonumber \\ &+& \sum^{\infty}_{k=1}
\frac{h_k^{(6)}}{k(2k-1)} - ln(2)\frac{h_k^{(5)}}{k(2k-1)} + ln(2)\sum^{\infty}_{k=1}\frac{H_k}{(2k-1)^5} + \frac{1}{2}\sum^{\infty}_{k=1} \frac{H_k}{k(2k-1)^6} + \frac{1}{2}
\sum^{\infty}_{k=1} \frac{H_k^2}{(2k-1)^5} \nonumber \\ &-& ln(2)\sum^{\infty}_{k=1}\frac{1}{k(2k-1)^5} - \frac{1}{2}\sum^{\infty}_{k=1}\frac{1}{k^2(2k-1)^6} -
\frac{1}{2}\sum^{\infty}_{k=1} \frac{h_k^{(5)}}{k^2(2k-1)^2} - 2\sum^{\infty}_{k=1} \frac{H_k h_k}{(2k-1)^5} \nonumber \\ &+& 2\sum^{\infty}_{k=1} \frac{h_k h_k^{(5)}}{k(2k-1)}~.
\end{eqnarray}
After elementary but tedious algebraic manipulations Eq.~(124) results. Also in this case one can observe that members of the fifth and sixth family are needed for an explicit calculation.

The second sum has been calculated to: 
\begin{eqnarray}
\sum^{\infty}_{k=1} \frac{H_k^{(2)} h_k^{(4)}}{k(2k-1)} &=& \frac{2921}{32}\zeta(7) + \frac{885}{32}ln(2)\zeta(6) - \frac{495}{32}\zeta(6) - \frac{1147}{32}\zeta(2)\zeta(5) + 
\frac{155}{32}\zeta(5) - \frac{51}{2}\zeta(3)\zeta(4) \nonumber \\ &-& \frac{15}{2}ln(2)\zeta(4) - \frac{63}{4}ln(2)\zeta(3)^2 + \frac{63}{8}\zeta(3)^2 -\frac{9}{2}\zeta(2)\zeta(3)
+ \frac{7}{2}\zeta(3)\sum^{\infty}_{k=1}\frac{h_k}{k^3} - 2ln(2)\sum^{\infty}_{k=1}\frac{h_k}{k^5} + \sum^{\infty}_{k=1}\frac{h_k}{k^5}  \nonumber \\ &-&
\frac{7}{2}\sum^{\infty}_{k=1} \frac{H_k^{(2)} h_k}{k^4} - 2\sum^{\infty}_{k=1} \frac{H_k h_k}{k^5}~.
\end{eqnarray}
In this case and for Eq.~(131) and also for Eq.~(134) \cite{ade16} is very helpful. Here the calculation starts with:
\begin{eqnarray}
\sum^{k}_{i=1} \frac{h_i^{(4)}}{i^2} = H_k^{(2)} h_k^{(4)}  - \sum^{k}_{i=1} \frac{H_{i-1}^{(2)}}{(2i-1)^4}~. 
\end{eqnarray}
From this it follows:
\begin{eqnarray}
\sum^{\infty}_{k=1} \frac{1}{k(2k-1)} \left( \sum^{k}_{i=1} \frac{h_i^{(4)}}{i^2} \right) &=& \sum^{\infty}_{k=1} \frac{H_k^{(2)} h_k^{(4)}}{k(2k-1)} -
\sum^{\infty}_{k=1} \frac{1}{k(2k-1)} \left( \sum^{k}_{i=1} \frac{H_{i-1}^{(2)}}{(2i-1)^4} \right)~.
\end{eqnarray}
Again, after lengthly algebraic manipulations Eq.~(128) results, where members of the first, second, fifth and seventh family are needed for an explicit calculation.

Furthermore we found:
\begin{eqnarray}
\sum^{\infty}_{k=1} \frac{H_k^{(3)} h_k^{(3)}}{k(2k-1)} &=& -\frac{2413}{32}\zeta(7) - \frac{135}{4}ln(2)\zeta(6) + \frac{135}{8}\zeta(6) + \frac{31}{2}\zeta(2)\zeta(5) - 
\frac{93}{2}\zeta(5) + \frac{489}{16}\zeta(3)\zeta(4) + \frac{75}{4}\zeta(4) \nonumber \\ &+& 21ln(2)\zeta(3)^2 - \frac{91}{8}\zeta(3)^2 +
21\zeta(2)\zeta(3) - 14ln(2)\zeta(3) - \frac{3}{2}\zeta(3)\sum^{\infty}_{k=1}\frac{h_k}{k^3} + 2\sum^{\infty}_{k=1} \frac{h_k}{k^3} \nonumber \\ &+& 
6ln(2)\sum^{\infty}_{k=1}\frac{h_k}{k^5} - 3\sum^{\infty}_{k=1}\frac{h_k}{k^5} + 3\sum^{\infty}_{k=1} \frac{H_k^{(2)} h_k}{k^4} +
6\sum^{\infty}_{k=1} \frac{H_k h_k}{k^5}~.
\end{eqnarray}
Here the calculation starts with:
\begin{eqnarray}
\sum^{k}_{i=1} \frac{h_i}{i^3} = H_k^{(3)} h_k  - \sum^{k}_{i=1} \frac{H_{i-1}^{(3)}}{(2i-1)}~. 
\end{eqnarray}
From this it follows:
\begin{eqnarray}
\sum^{\infty}_{k=1} \frac{1}{(2k-1)^3} \left( \sum^{k}_{i=1} \frac{h_i}{i^3} \right) &=& \sum^{\infty}_{k=1} \frac{H_k^{(3)} h_k}{(2k-1)^3} -
\sum^{\infty}_{k=1} \frac{1}{(2k-1)^3} \left( \sum^{k}_{i=1} \frac{H_{i-1}^{(3)}}{(2i-1)} \right)~.
\end{eqnarray}
Again, after lengthly algebraic manipulations Eq.~(131) results, where members of the first, second, fifth and seventh family are needed for an explicit calculation.

Next to last it follows:
\begin{eqnarray}
\sum^{\infty}_{k=1} \frac{H_k^{(4)} h_k^{(2)}}{k(2k-1)} &=& \frac{381}{4}\zeta(7) + \frac{159}{8}ln(2)\zeta(6) - \frac{45}{4}\zeta(6) - \frac{31}{8}\zeta(2)\zeta(5) +
62\zeta(5) - \frac{223}{4}\zeta(3)\zeta(4) - \frac{45}{2}\zeta(4) \nonumber \\ &-& 7ln(2)\zeta(3)^2 + \frac{7}{2}\zeta(3)^2 - 35\zeta(2)\zeta(3) + 42\zeta(3) - 24ln(2)\zeta(2) -
2\zeta(3)\sum^{\infty}_{k=1}\frac{h_k}{k^3} + 8\sum^{\infty}_{k=1} \frac{h_k}{k^3} \nonumber \\ &-& 8ln(2)\sum^{\infty}_{k=1}\frac{h_k}{k^5} + 4\sum^{\infty}_{k=1}\frac{h_k}{k^5} -
6\sum^{\infty}_{k=1} \frac{H_k^{(2)} h_k}{k^4} - 8\sum^{\infty}_{k=1} \frac{H_k h_k}{k^5}~.
\end{eqnarray}
Here the calculation starts with:
\begin{eqnarray}
\sum^{k}_{i=1} \frac{h_i^{(2)}}{i^4} = H_k^{(2)} h_k^{(4)} - \sum^{k}_{i=1} \frac{H_{i-1}^{(2)}}{(2i-1)^{(4)}}~. 
\end{eqnarray}
From this it follows:
\begin{eqnarray}
\sum^{\infty}_{k=1} \frac{1}{k(2k-1)} \left( \sum^{k}_{i=1} \frac{h_i^{(2)}}{i^4} \right) &=& \sum^{\infty}_{k=1} \frac{H_k^{(4)} h_k^{(2)}}{k(2k-1)} -
\sum^{\infty}_{k=1} \frac{1}{k(2k-1)} \left( \sum^{k}_{i=1} \frac{H_{i-1}^{(4)}}{(2i-1)^2} \right)~.
\end{eqnarray}
Again, after lengthly algebraic manipulations Eq.~(134) results, where members of the first, second, fifth and seventh family are needed for an explicit calculation.

Finally the last sum of this family results to:
\begin{eqnarray}
\sum^{\infty}_{k=1} \frac{H_k^{(5)} h_k}{k(2k-1)} &=& -\frac{127}{64}\zeta(7) - \frac{209}{8}\zeta(2)\zeta(5) - 31\zeta(5) + \frac{455}{16}\zeta(3)\zeta(4) + 14\zeta(2)\zeta(3)
-28\zeta(3) + 32\zeta(2) + \zeta(3)\sum^{\infty}_{k=1}\frac{h_k}{k^3} \nonumber \\ &-& 8\sum^{\infty}_{k=1}\frac{h_k}{k^3} + 2ln(2)\sum^{\infty}_{k=1}\frac{h_k}{k^5} -
2\sum^{\infty}_{k=1}\frac{h_k}{k^5} + \frac{7}{2}\sum^{\infty}_{k=1} \frac{H_k^{(2)} h_k}{k^4} + 2\sum^{\infty}_{k=1} \frac{H_k h_k}{k^5}~.
\end{eqnarray}
Again, in this case \cite{ade16} is very helpful. We start the calculation with:
\begin{eqnarray}
H_k^{(5)} h_k = \sum^{k}_{i=1} \frac{h_i}{i^5} + \sum^{k}_{i=1}\frac{H_{i-1}^{(5)}}{(2i-1)}~.
\end{eqnarray} 
From this we get:
\begin{eqnarray}
\sum^{\infty}_{k=1} \frac{H_k^{(5)}h_k}{k(2k-1)} = \sum^{\infty}_{k=1} \frac{1}{k(2k-1)} \left( \sum^{k}_{i=1} \frac{h_i}{i^5} \right) + \sum^{\infty}_{k=1}\frac{1}{k(2k-1)} 
\left( \frac{H_{i-1}^{(5)}}{(2i-1)} \right)~.
\end{eqnarray}
Furthermore it follows:
\begin{eqnarray}
\sum^{\infty}_{k=1} \frac{H_k^{(5)} h_k}{k(2k-1)} &=& ln(2)\sum^{\infty}_{k=1} \frac{h_k}{k^5} - \frac{1}{2}ln(2)\sum^{\infty}_{k=1}\frac{H_k}{k^5} - \frac{1}{4}\sum^{\infty}_{k=1}\frac{H_k^2}{k^5}
- \frac{1}{2}\sum^{\infty}_{k=1}\frac{H_k^{(5)}}{k(2k-1)} + \frac{1}{4}\sum^{\infty}_{k=1}\frac{H_k^{(5)}}{k^2(2k-1)^2} - \frac{1}{2}\sum^{\infty}_{k=1}\frac{H_k}{k(2k-1)^5} 
\nonumber \\ &+& \frac{1}{4}\sum^{\infty}_{k=1}\frac{H_k^{(6)}}{k(2k-1)} - \frac{1}{4}\sum^{\infty}_{k=1}\frac{H_k}{k^6(2k-1)} + \sum^{\infty}_{k=1}\frac{h_k}{k(2k-1)^5} + \frac{1}{2}
\sum^{\infty}_{k=1}\frac{h_k}{k^6(2k-1)} - ln(2)\sum^{\infty}_{k=1}\frac{1}{k(2k-1)^5} \nonumber \\ &-& \frac{1}{2}\sum^{\infty}_{k=1}\frac{h_k}{k^2(2k-1)^6} -
\sum^{\infty}_{k=1}\frac{h_k^2}{k^5} + \sum^{\infty}_{k=1}\frac{H_k h_k}{k^5} - \sum^{\infty}_{k=1} \frac{H_k H_k^{(5)}}{k(2k-1)}~.
\end{eqnarray}
Finally, after lengthly algebraic manipulations Eq.~(137) results, where members of the first, third, and seventh family are needed for an explicit calculation.
With this we have calculated all members of the seventh family.

\subsection{Eighth family}
Here we present the following results:
\begin{eqnarray}
\sum^{\infty}_{k=1} \frac{h_k^{2}}{k^5} &=& \frac{4191}{64}\zeta(7) - \frac{155}{8}\zeta(2)\zeta(5) - \frac{343}{16}\zeta(3)\zeta(4) - 
\frac{7}{2}\sum^{\infty}_{k=1} \frac{H_k^{(2)} h_k}{k^4}~,
\end{eqnarray}
In order to prove this result, we introduce the following two valued help function:
\begin{eqnarray}
\sum^{k}_{i=1}\frac{h_i}{i+k} = \sum^{k}_{i=1}\frac{h_i}{i} - \frac{1}{2}H_k h_k~. 
\end{eqnarray}
With this one gets 
\begin{eqnarray}
\sum^{\infty}_{k=1} \frac{1}{k^5}\left(\sum^{k}_{i=1}\frac{h_i}{i+k} \right) = \sum^{\infty}_{k=1} \frac{1}{k^5}\left(\sum^{k}_{i=1}\frac{h_i}{i} \right) - \frac{1}{2}
\sum^{\infty}_{k=1} \frac{H_k h_k}{k^5}~.
\end{eqnarray}
Rearranging the summation order in the sum of the left side of Eq.~(143) it follows:
\begin{eqnarray}
\sum^{\infty}_{k=1} \frac{1}{k^5}\left(\sum^{k}_{i=1}\frac{h_i}{i+k} \right) = \frac{1}{2}\sum^{\infty}_{k=1} \frac{h_k}{k^6} + \sum^{\infty}_{k=1} h_k \left( \sum^{\infty}_{i=1} 
\frac{1}{(i+k)^3(i+2k)} \right)~.
\end{eqnarray}
Applying a partial fraction decomposition the result is:
\begin{eqnarray}
\sum^{\infty}_{k=1} \frac{1}{k^5}\left(\sum^{k}_{i=1}\frac{h_i}{i+k} \right) &=& \frac{1}{2}\sum^{\infty}_{k=1} \frac{h_k}{k^6} + \sum^{\infty}_{k=1} \frac{h_k}{k}
\left(\zeta(5) - H_k^{(5)} \right) - \sum^{\infty}_{k=1} \frac{h_k}{k^2} \left(\zeta(4) - H_k^{(4)} \right) + \sum^{\infty}_{k=1} \frac{h_k}{k^3} \left(\zeta(3) - H_k^{(3)} \right)
\nonumber \\ &-& \sum^{\infty}_{k=1} \frac{h_k}{k^4} \left(\zeta(2) - H_k^{(2)} \right) + \sum^{\infty}_{k=1} \frac{h_k}{k^4} \left(\frac{h_k}{k} - \frac{1}{2} \frac{H_k}{k} \right)
\end{eqnarray}
Explicit calculation and comparison with Eq.~(143) results to:
\begin{eqnarray}
\sum^{\infty}_{k=1} \frac{h_k^{2}}{k^5} &=& \frac{1}{2}\sum^{\infty}_{k=1} \frac{h_k}{k^6} + \zeta(4)\frac{h_k}{k^2} - \zeta(3)\frac{h_k}{k^3} + \zeta(2)\frac{h_k}{k^4} -
\sum^{\infty}_{k=1} \frac{H_k^{(2)} h_k}{k^4} - \sum^{\infty}_{k=1} \frac{H_k^{(4)} h_k}{k^2} + \sum^{\infty}_{k=1} \frac{H_k^{(3)} h_k}{k^3}~.
\end{eqnarray}
Using Eq.~(42) and (43) we arrive at the final result, namely Eq.~(141)

The next identity is:
\begin{eqnarray}
\sum^{\infty}_{k=1} \frac{h_k h_k^{(2)}}{k^4} &=& -\frac{6477}{128}\zeta(7) + \frac{155}{8}\zeta(2)\zeta(5) + \frac{371}{32}\zeta(3)\zeta(4) + 
\frac{7}{2}\sum^{\infty}_{k=1} \frac{H_k^{(2)} h_k}{k^4}~.
\end{eqnarray}
This Euler sum is closely connected with the Euler sum presented in Eq.~(163). We start the corresponding calculation with:
\begin{eqnarray}
\sum^{\infty}_{k=1}\frac{1}{(2k-1)^2} \left( \sum^{\infty}_{i=1} \frac{h_i^{(4)}}{(2k-1)(2i+2k-1)} \right) &=& \frac{3}{8}\zeta(2)\sum^{\infty}_{k=1}\frac{h_k^{(4)}}{k(2k-1)}
-\frac{1}{2}\sum^{\infty}_{k=1}\frac{h_k^{(4)}}{k(2k-1)} \left( \sum^{\infty}_{i=1} \frac{1}{(2k-1)(2i+2k-1)} \right) \nonumber \\  &=&
\frac{3}{8}\zeta(2)\sum^{\infty}_{k=1}\frac{h_k^{(4)}}{k(2k-1)} - \frac{1}{4}\sum^{\infty}_{k=1}\frac{h_k h_k^{(4)}}{k^2 (2k-1)}~. 
\end{eqnarray}
Using Lemma 4a from \cite{bra21} it follows:
\begin{eqnarray}
&&\frac{31}{64}\zeta(5)\sum^{\infty}_{k=1}\frac{1}{k(2k-1)^2} + \frac{15}{64}\zeta(4)\sum^{\infty}_{k=1}\frac{H_{k-1}}{k(2k-1)^2} - \frac{7}{64}\zeta(3)\sum^{\infty}_{k=1}
\frac{H^{(2)}_{k-1}}{k(2k-1)^2} + \frac{3}{64}\zeta(2)\sum^{\infty}_{k=1} \frac{H^{(3)}_{k-1}}{k(2k-1)^2} \nonumber \\ &-& \frac{3}{64}\zeta(2)\sum^{\infty}_{k=1}\frac{h_k}{k^4} +
\frac{1}{16} \sum^{\infty}_{k=1}\frac{h_k h_k^{(2)}}{k^4} - \frac{1}{16} \sum^{\infty}_{k=1}\frac{h_k^2}{k^4} + \frac{1}{32} \sum^{\infty}_{k=1}\frac{H_k h_k}{k^4} =
\frac{3}{8}\zeta(2)\sum^{\infty}_{k=1}\frac{h_k^{(4)}}{k(2k-1)} - \frac{1}{2}\zeta(2)\sum^{\infty}_{k=1}\frac{h_k h_k^{(4)}}{k(2k-1)}  \nonumber \\ &+&
\frac{1}{4}\sum^{\infty}_{k=1}\frac{h_k h_k^{(4)}}{k^2}~.
\end{eqnarray}
After a very lengthy computation we finally arrive at:
\begin{eqnarray}
\sum^{\infty}_{k=1}\frac{h_k h_k^{(4)}}{k^2} = \frac{315}{128}\zeta(3)\zeta(4) - \frac{7}{8}\zeta(3)\sum^{\infty}_{k=1}\frac{h_k}{k^3} +
\frac{1}{4} \sum^{\infty}_{k=1}\frac{h_k h_k^{(2)}}{k^4}~. 
\end{eqnarray}
As a consequence one is able to compute the Euler sum in Eq.~(163) if Eq.~(147) is known or vice versa. This example may shown again that the ranking in the different calculations
is not arbitrary but highly recursive where different Euler sums from different families are involved.

Furthermore we have:
\begin{eqnarray}
\sum^{\infty}_{k=1} \frac{h_k h_k^{(3)}}{k^3} &=& -\frac{1905}{128}\zeta(7) + \frac{93}{16}\zeta(2)\zeta(5) + \frac{315}{64}\zeta(3)\zeta(4)~,
\end{eqnarray}
As a third example we present the calculational procedure that leads to Eq.~(151). First we introduce the two valued help function
\begin{eqnarray}
\sum^{k}_{i=1}\frac{h_i^{(3)}}{i+k} =  \frac{1}{2}H_k h_k^{(3)} + h_k^{(2)}h_k^{(2)} - h_k h_k^{(3)} - \sum^{k}_{i=1} \frac{H_{i-1}}{(2i-1)^3}~. 
\end{eqnarray}
With this we get:
\begin{eqnarray}
\sum^{\infty}_{k=1} \frac{1}{k^3}\left( \sum^{k}_{i=1} \frac{h_i^{(3)}}{i+k} \right) &=& \frac{1}{2} \sum^{\infty}_{k=1} \frac{H_k h_k^{(3)}}{k^3} +
\sum^{\infty}_{k=1} \frac{h_k^{(2)} h_k^{(2)}}{k^3} - \sum^{\infty}_{k=1} \frac{h_k h_k^{(3)}}{k^3} - \sum^{\infty}_{k=1} \frac{1}{k^3}
\left( \sum^{k}_{i=1} \frac{H_{i-1}}{(2i-1)^3}  \right)~.
\end{eqnarray}
Rearanging the summation order on the left side it follows:
\begin{eqnarray}
\sum^{\infty}_{k=1} \frac{1}{k^3}\left( \sum^{k}_{i=1} \frac{h_i^{(3)}}{i+k} \right) &=& \frac{1}{2} \sum^{\infty}_{k=1} \frac{h_k^{(3)}}{k^4} +
\sum^{\infty}_{k=1} h_k^{(3)} \left( \sum^{\infty}_{i=1} \frac{1}{i(i+k)^3(i+2k)} \right)~.
\end{eqnarray}
From this it follows:
\begin{eqnarray}
\sum^{\infty}_{k=1} \frac{1}{k^3}\left( \sum^{k}_{i=1} \frac{h_i^{(3)}}{i+k} \right) &=& \frac{1}{2} \sum^{\infty}_{k=1} \frac{h_k^{(3)}}{k^4} +
\sum^{\infty}_{k=1} \frac{h_k^{(3)}}{k}\left(\zeta(3)-H_k^{(3)}\right) - \sum^{\infty}_{k=1} \frac{h_k^{(3)}}{k^2}\left(\zeta(2)-H_k^{(2)} \right) \nonumber \\ &+&
\sum^{\infty}_{k=1} \frac{h_k^{(3)}}{k^2}\left( \frac{H_{2k}}{k}-\frac{H_{k}}{k} \right)~.
\end{eqnarray}
Further computation results in:
\begin{eqnarray}
2\sum^{\infty}_{k=1} \frac{h_k h_k^{(3)}}{k^3} &=& \frac{1}{2} \sum^{\infty}_{k=1} \frac{h_k^{(3)}}{k^4} + \zeta(2)\sum^{\infty}_{k=1} \frac{h_k^{(3)}}{k^2} + 
\sum^{\infty}_{k=1} \frac{H_k h_k^{(3)}}{k^3} - \sum^{\infty}_{k=1} \frac{H_k^{(2)} h_k^{(3)}}{k^2} - \sum^{\infty}_{k=1} \frac{h_k^{(2)} h_k^{(2)}}{k^3}
- \sum^{\infty}_{k=1} \frac{1}{k^3} \left( \sum^{k}_{i=1} \frac{H_{i-1}}{(2i-1)^3}  \right) \nonumber \\ &-& \sum^{\infty}_{k=1} \frac{1}{k^3}
\left( \sum^{k}_{i=1} \frac{h_i^{(3)}}{i} \right)~, 
\end{eqnarray}
and finally it follows:
\begin{eqnarray}
\sum^{\infty}_{k=1} \frac{h_k h_k^{(3)}}{k^3} &=& \frac{1}{4} \sum^{\infty}_{k=1} \frac{h_k^{(3)}}{k^4} + \frac{1}{2} \zeta(2)\sum^{\infty}_{k=1} \frac{h_k^{(3)}}{k^2} -
\frac{1}{2} \sum^{\infty}_{k=1} \frac{H_k^{(2)} h_k^{(3)}}{k^2} + \frac{1}{2} \sum^{\infty}_{k=1} \frac{h_k^{(2)} h_k^{(2)}}{k^3}~. 
\end{eqnarray}
Using Eq.~(49) and Eq.~(158) the expression presented in Eq.~(151) follows. This computational procedure demonstrates that two other Euler sums belonging to the first and eighth
family must be known in order to be able to calculate explicitly Eq.~(151).

Also we found:
\begin{eqnarray}
\sum^{\infty}_{k=1} \frac{h_k^{(2)} h_k^{(2)}}{k^3} &=& \frac{5207}{128}\zeta(7) - \frac{217}{8}\zeta(2)\zeta(5) + \frac{259}{32}\zeta(3)\zeta(4) - 
\frac{7}{2}\sum^{\infty}_{k=1} \frac{H_k^{(2)} h_k}{k^4}~.
\end{eqnarray}
and
\begin{eqnarray}
\sum^{\infty}_{k=1} \frac{h_k^{(2)} h_k^{(3)}}{k^2} &=& \frac{16637}{512}\zeta(7) - \frac{217}{32}\zeta(2)\zeta(5) - \frac{973}{64}\zeta(3)\zeta(4) + 
\frac{7}{8}\zeta(3)\sum^{\infty}_{k=1}\frac{h_k}{k^3} - \frac{7}{8}\sum^{\infty}_{k=1} \frac{H_k^{(2)} h_k}{k^4}~.
\end{eqnarray}
The last computation concerns Eq.~(159). Here we start with the following two valued help function:
\begin{eqnarray}
\sum^{\infty}_{i=1} \frac{h_i^{(3)}}{(2k+2i-1)^2} = \sum^{\infty}_{i=1} \frac{h_k^{(3)}}{(2k+1)^2} + \frac{9}{32}\zeta(2)H_{k-1}^{(3)} - \frac{7}{32}\zeta(2)H_{k-1}^{(2)}
- \frac{3}{16} \sum^{k-1}_{i=1}\frac{h_i}{i^4} - \frac{1}{8}\sum^{k-1}_{i=1}\frac{h_i^{(2)}}{i^3}~.
\end{eqnarray}
With this we get:
\begin{eqnarray}
\sum^{\infty}_{k=1}\frac{1}{(2k-1)^2} \left( \sum^{\infty}_{i=1} \frac{h_i^{(3)}}{(2k+2i-1)^2} \right) &=& \sum^{\infty}_{i=1} h_i^{(3)} \left( 
\sum^{\infty}_{k=1}\frac{1}{(2k-1)^2 (2k+2i-1)^2} \right) \nonumber \\  &=& \frac{3}{8}\zeta(2)\sum^{\infty}_{k=1}\frac{h_k^{(3)}}{k^2} - \frac{1}{4}\sum^{\infty}_{k=1}\frac{h_k h_k^{(3)}}{k^3}
- \frac{1}{4}\sum^{\infty}_{k=1}\frac{h_k^{(2)} h_k^{(3)}}{k^2}~.
\end{eqnarray}
Calculating explicitly the left side of Eq.~(161) it follows: 
\begin{eqnarray}
&&\frac{3}{8}\zeta(2)\sum^{\infty}_{k=1}\frac{h_k^{(3)}}{k^2} - \frac{1}{4}\sum^{\infty}_{k=1}\frac{h_k h_k^{(3)}}{k^3} - \frac{1}{4}\sum^{\infty}_{k=1}\frac{h_k^{(2)} h_k^{(3)}}{k^2} =
\frac{3}{4}\zeta(2)\sum^{\infty}_{k=1} \frac{h_k^{(3)}}{(2k+1)^2} + \frac{9}{32}\zeta(2)\sum^{\infty}_{k=1}\frac{H_{k-1}^{(3)}}{(2k-1)^2} -
\nonumber \\ && \frac{7}{32}\zeta(3)\sum^{\infty}_{k=1}\frac{H_{k-1}^{(2)}}{(2k-1)^2} - \frac{9}{64}\zeta(2)\sum^{\infty}_{k=1} \frac{h_k}{k^4} -
\frac{3}{32}\zeta(2)\frac{h_k^{(2)}}{k^3} + \frac{3}{16}\sum^{\infty}_{k=1} \frac{h_k h_k^{(2)}}{k^4} +\frac{1}{8} \sum^{\infty}_{k=1} \frac{h_k^{(2)} h_k^{(2)}}{k^3}~.
\end{eqnarray}
Inspecting Eq.~(162) one can observe that for the explicit calculation of Eq.~(159), Eq.~(147), Eq.~(151) and Eq.~(158) are required. Again a highly recursive situation.

Finally it follows:
\begin{eqnarray}
\sum^{\infty}_{k=1} \frac{h_k h_k^{(4)}}{k^2} &=& -\frac{6477}{512}\zeta(7) + \frac{155}{32}\zeta(2)\zeta(5) + \frac{343}{64}\zeta(3)\zeta(4) - 
\frac{7}{8}\zeta(3)\sum^{\infty}_{k=1}\frac{h_k}{k^3} + \frac{7}{8}\sum^{\infty}_{k=1} \frac{H_k^{(2)} h_k}{k^4}~.
\end{eqnarray}
With this we have calculated all members of the eighth family concerning the order seven.

\section{Summary}
We have introduced a special summation method that allows to calculate explicitly eighth families of quadratic Euler sums in terms of zeta values and polylogarithmic values $Li_n(1/2)$.
Our approach is able to calculate explicitly the corresponding families of nonlinear Euler sums of odd orders. The order seven case seems to be more complicated as higher order polylogarithmic
values like $Li_7(1/2)$ are not known explicitly but most probably are necessary for the complete calculation of the nonlinear Euler sums being members of the first family in the order seven or
even higher order cases. For a complete reduction procedure in the order seven case it remains to find a direct relation between $\sum^{\infty}_{k=1} \frac{H_k^{(2)} h_k}{k^4}$
and $\sum^{\infty}_{k=1} \frac{H_k h_k}{k^5}$ in terms of zeta values and linear Euler sums. Up to now only the following approximation to $\sum^{\infty}_{k=1} \frac{H_k^{(2)} h_k}{k^4}$ is
known:
\begin{eqnarray}
\sum^{\infty}_{k=1} \frac{H_k^{(2)} h_k}{k^4} \approx -\frac{1559}{1943}\zeta(2)\zeta(5) + \frac{1469}{759}\zeta(3)\zeta(4)~.
\end{eqnarray}
Remarkably, this results is correct within $2.9*10^{-15}$.
If one inspects Eq.~(15) it seems plausible that an explicit calculation of $\sum^{\infty}_{k=1}
\frac{H_k h_k}{k^5}$ in terms of zeta values, linear Euler sums and $Li_7(1/2)$ should be possible. This would completes the reduction procedure which can be applied to the order nine or even
higher odd order cases if $Li_n(1/2)$ with n=7,9, ... is known. Nevertheless a significant reduction of quadratic Euler sums to zeta values and linear Euler sums is possible.

\section{Appendix A}
Here we present as an outlook few related ternary Euler sums of different orders. In a forthcoming paper we will show that even for a variety of ternary Euler sums belonging to the different families a
systematic reduction process is possible.

First we present some examples for the order 4 case:
\begin{eqnarray}
\sum^{\infty}_{k=1} \frac{H_k^3}{k(2k-1)} &=& -\frac{37}{8}\zeta(4) -7ln(2)\zeta(3) + 7\zeta(3) + 6\left( ln(2) \right)^2\zeta(2) - 12ln(2)\zeta(2) + 8\zeta(2)
+ 16ln(2)- 24\left( ln(2) \right)^2 \nonumber \\ &+& 16\left( ln(2) \right)^3-4\left( ln(2) \right)^4 + \sum^{\infty}_{k=1} \frac{h_k}{k^3}
\end{eqnarray}

\begin{eqnarray}
\sum^{\infty}_{k=1} \frac{h_k^3}{k(2k-1)} &=& \frac{57}{64}\zeta(4) + \frac{7}{8}ln(2)\zeta(3) + \frac{9}{4}\left(ln(2)\right)^2\zeta(2) 
- \frac{1}{8}\sum^{\infty}_{k=1}\frac{h_k}{k^3}~.
\end{eqnarray}

and
\begin{eqnarray}
\sum^{\infty}_{k=1} \frac{H_k^2 h_k}{k(2k-1)} &=& \frac{1}{4}\zeta(4) - \frac{17}{2}ln(2)\zeta(3) + \frac{17}{2}\zeta(3) + 4\left(ln(2)\right)^2\zeta(2) 
- 8ln(2)\zeta(2) + 4\zeta(2) + \sum^{\infty}_{k=1}\frac{h_k}{k^3}~.
\end{eqnarray}

Next we present two examples for the order 5 case:
\begin{eqnarray}
\sum^{\infty}_{k=1} \frac{H_k h_k^2}{k^2} &=& \frac{155}{32}\zeta(5) + \frac{7}{8}\zeta(2)\zeta(3)
\end{eqnarray}

and
\begin{eqnarray}
\sum^{\infty}_{k=1} \frac{h_k^3}{(2k-1)^2} &=& -\frac{713}{256}\zeta(5) + \frac{249}{64}ln(2)\zeta(4) + \frac{3}{32}\zeta(2)\zeta(3) + \left(ln(2)\right)^3\zeta(2) 
- \frac{1}{40}\left(ln(2)\right)^5 + 3Li_5\left(\frac{1}{2}\right)~.
\end{eqnarray}

For the order 6 case we found:
\begin{eqnarray}
\sum^{\infty}_{k=1} \frac{H_k^3}{(2k-1)^3} &=& \frac{405}{128}\zeta(6) - \frac{93}{4}ln(2)\zeta(5) + \frac{93}{4}\zeta(5) + \frac{135}{8}\left(ln(2)\right)^2\zeta(4) - 
\frac{135}{4}ln(2)\zeta(4) - \frac{7}{32}\zeta(3)^2 + \frac{21}{4}ln(2)\zeta(2)\zeta(3) \nonumber \\ &-& \frac{21}{4}\zeta(2)\zeta(3) - 7\left(ln(2)\right)^3\zeta(3)
+ 21\left(ln(2)\right)^2\zeta(3) + 21ln(2)\zeta(3) - \frac{63}{2}\zeta(3) - 18\left(ln(2)\right)^2\zeta(2) \nonumber \\ &+& 30ln(2)\zeta(2) - 6\zeta(2) + 48ln(2)-   
36\left(ln(2)\right)^2 + 8\left(ln(2)\right)^3 - 6\sum^{\infty}_{k=1}\frac{h_k}{k^3} + \frac{3}{4}\sum^{\infty}_{k=1}\frac{h_k}{k^5}~.
\end{eqnarray}

For the order 7 case we present the following sums:
\begin{eqnarray}
\sum^{\infty}_{k=1} \frac{h_k^3}{k^4} =  \frac{93}{4}\zeta(2)\zeta(5) - \frac{945}{32}\zeta(3)\zeta(4)
\end{eqnarray}

\begin{eqnarray}
\sum^{\infty}_{k=1} \frac{h_k^{(2)} h_k^2}{k^3} =  -\frac{93}{8}\zeta(2)\zeta(5) + \frac{525}{32}\zeta(3)\zeta(4)
\end{eqnarray}

\begin{eqnarray}
\sum^{\infty}_{k=1} \frac{h_k^3}{(2k-1)^4} &=& \frac{127}{256}\zeta(7) + \frac{441}{512}ln(2)\zeta(6) + \frac{3}{64}\zeta(2)\zeta(5) + \frac{93}{64}\left(ln(2)\right)^2\zeta(5) 
- \frac{195}{512}\zeta(3)\zeta(4) + \frac{15}{16}\left(ln(2)\right)^3\zeta(4) \nonumber \\ &-& \frac{21}{64}ln(2)\zeta(3)^2 - \frac{9}{32}\left(ln(2)\right)^2\zeta(2)\zeta(3) 
- \frac{3}{32}ln(2)\sum^{\infty}_{k=1}\frac{h_k}{k^5} - \frac{3}{64}\sum^{\infty}_{k=1} \frac{H_k h_k}{k^5}~,
\end{eqnarray}

\begin{eqnarray}
\sum^{\infty}_{k=1} \frac{h_k h_k^{(2)} h_k^{(2)}}{k^2} &=& \frac{1113}{64}ln(2)\zeta(6) - \frac{1023}{128}\zeta(2)\zeta(5) - \frac{105}{32}\zeta(3)\zeta(4) + 
\frac{5}{2}\left(ln(2)\right)^3\zeta(4) - \frac{21}{4}\left(ln(2)\right)^2\zeta(2)\zeta(3) \nonumber \\ &-& \frac{1}{10}\left(ln(2)\right)^5\zeta(2) +
\frac{3}{2}ln(2)\zeta(2)\sum^{\infty}_{k=1}\frac{h_k}{k^3} + 12\zeta(2) Li_5\left(\frac{1}{2}\right)~, 
\end{eqnarray}

\begin{eqnarray}
\sum^{\infty}_{k=1} \frac{h_k^{(3)} h_k^2}{k^2} &=& -\frac{1113}{64}ln(2)\zeta(6) + \frac{1767}{128}\zeta(2)\zeta(5) - \frac{5}{2}\left(ln(2)\right)^3\zeta(4) +
\frac{21}{4}\left(ln(2)\right)^2\zeta(2)\zeta(3) + \frac{1}{10}\left(ln(2)\right)^5\zeta(2) \nonumber \\ &-& \frac{3}{2}ln(2)\zeta(2)\sum^{\infty}_{k=1}\frac{h_k}{k^3} - 
12\zeta(2) Li_5\left(\frac{1}{2}\right)~, 
\end{eqnarray}

\begin{eqnarray}
\sum^{\infty}_{k=1} \frac{H_k^{(2)} h_k^2}{k^3} = \frac{15113}{128}\zeta(7) - \frac{155}{4}\zeta(2)\zeta(5) - \frac{1249}{32}\zeta(3)\zeta(4) + \frac{7}{2}\zeta(3)
\sum^{\infty}_{k=1}\frac{h_k}{k^3} - \frac{21}{4}\sum^{\infty}_{k=1} \frac{H_k^{(2)} h_k}{k^4}~,
\end{eqnarray}

\section{Appendix B}
Finally we present two biquadratic Euler sums belonging to the third and second family:
\begin{eqnarray}
\sum^{\infty}_{k=1} \frac{H_k^4}{k(2k-1)} &=& -\frac{305}{4}\zeta(5) + 47ln(2)\zeta(4) + 6\zeta(4) - 3\zeta(2)\zeta(3) - 56ln(2)\zeta(3) + 40\zeta(3) -\frac{32}{3}\left(ln(2)\right)^3\zeta(2)
\nonumber \\ &+& 48\left(ln(2)\right)^3\zeta(2) - 48ln(2)\zeta(2) + 24\zeta(2) + 32ln(2) - 64\left(ln(2)\right)^2 +64\left(ln(2)\right)^3 - 32\left(ln(2)\right)^4  \nonumber \\ &+&
\frac{88}{15}\left(ln(2)\right)^5 + 8\sum^{\infty}_{k=1}\frac{h_k}{k^3} +64Li_5\left(\frac{1}{2}\right)~.
\end{eqnarray}
and
\begin{eqnarray}
\sum^{\infty}_{k=1} \frac{H_k^4}{(2k-1)^3} &=& \frac{13081}{128}\zeta(7) - \frac{405}{16}ln(2)\zeta(6) + \frac{405}{16}\zeta(6)- \frac{403}{16}\zeta(2)\zeta(5) +
93\left(ln(2)\right)^2\zeta(5) - 186ln(2)\zeta(5) \nonumber \\ &+& \frac{1643}{4}\zeta(5) - \frac{1141}{32}\zeta(3)\zeta(4) - 45\left(ln(2)\right)^3\zeta(4) +
135\left(ln(2)\right)^2\zeta(4) - 318ln(2)\zeta(4) - 87\zeta(4) \nonumber \\ &+& \frac{7}{4}ln(2)\zeta(3)^2 - \frac{7}{4}\zeta(3)^2 - 21\left(ln(2)\right)^2 \zeta(2)\zeta(3)
+ 42ln(2)\zeta(2)\zeta(3) - 35\zeta(2)\zeta(3) + 14\left(ln(2)\right)^4\zeta(3)  \nonumber \\ &-&  56\left(ln(2)\right)^3\zeta(3) + 84\left(ln(2)\right)^2\zeta(3) +
252ln(2)\zeta(3) - 94\zeta(3) + 16\left(ln(2)\right)^3\zeta(2) - 120\left(ln(2)\right)^2\zeta(2) \nonumber \\ &+& 72ln(2)\zeta(2)+ 24\zeta(2) + 160ln(2) -   
192\left(ln(2)\right)^2 + 96\left(ln(2)\right)^3 -16\left(ln(2)\right)^4 + \frac{16}{5}\left(ln(2)\right)^5 \nonumber \\ &-& 44\sum^{\infty}_{k=1}\frac{h_k}{k^3} -
6ln(2)\sum^{\infty}_{k=1}\frac{h_k}{k^5} + 6\sum^{\infty}_{k=1}\frac{h_k}{k^5} - 3\sum^{\infty}_{k=1} \frac{H_k h_k}{k^5} - \frac{13}{4} \sum^{\infty}_{k=1} \frac{H_k^{(2)} h_k}{k^4}
- 384Li_5\left(\frac{1}{2}\right)~.
\end{eqnarray}

\end{document}